\documentclass[nofootinbib]{revtex4-2}  
\usepackage{graphicx}  
\usepackage{amsmath}
\usepackage{amssymb}
\usepackage{dcolumn}   
\usepackage{bm}        
\usepackage{amssymb}   
\usepackage{color}
\allowdisplaybreaks
\newtheorem{theorem}{Theorem}
\newtheorem{example}{Example}

\begin{document}
	
	\title{Adaptive Multiquadratic Radial Basis Function-based Explicit Runge--Kutta Methods} 
	\author{Rajesh Yadav}%
	\email{rajeshy@rgipt.ac.in}
	\author{Deepak Kumar Yadav}
	\email{20bs0006@rgipt.ac.in }
\affiliation{Department of Mathematical Sciences, Rajiv Gandhi Institute of Petroleum Technology, Jais, Amethi, 229304, India}
\author{Alpesh Kumar}%
\email{alpesh@allduniv.ac.in}
\affiliation{ Department of Mathematics, University of Allahabad, Prayagraj, Uttar Pradesh, 211002, India}

\begin{abstract}
\noindent 
Runge--Kutta (RK) methods are widely used techniques for solving a class of initial value problems. In this article, we introduce an adaptive multiquadratic (MQ) radial basis function (RBF)-based method to develop enhanced explicit RK methods. These methods achieve a higher order of convergence than the corresponding classical RK methods. To improve the local convergence of the numerical solution, we optimize the free parameters (shape functions) involved in the RBFs by forcing the local truncation errors to vanish. We also present a convergence and stability analysis of the proposed methods. To demonstrate the advantages of these methods in terms of accuracy and convergence, we consider several numerical examples and compare the performance of our methods with that of the classical RK methods. The Tables and Figures presented in this article clearly validate the superiority of the proposed methods.  
\end{abstract}
\maketitle
\section{Introduction}
\noindent
The numerical solution of initial value problems (IVPs) in ordinary differential equations (ODEs) is fundamental to computational mathematics, with wide-ranging applications in physics, engineering, biology, and control theory. Among the various numerical techniques, Runge–Kutta (RK) methods are particularly popular due to their self-starting nature, ease of implementation, and ability to achieve high accuracy without requiring excessively small step sizes. Since their introduction, numerous studies have focused on the generalization of classical RK methods for specific applications. Some notable recent examples include \cite{Hochbruck, Simos, Boscarino, Sebastiano,Bresten, Mei, Gu25}. 
The core idea behind RK methods is to compute a sequence of intermediate values within each time step to improve the accuracy of the approximation. For explicit RK methods, the number of stages $(s)$ typically determines the attainable order of convergence when $s\le 4$. In particular,  an explicit $s$ stage RK method can achieve $s$-th order accuracy for $s\le 4$. However, due to the constraints imposed by the Butcher order conditions, to achieve $s$-th order accuracy for $s>4$, the method generally requires more than $s$ stages. For example, a fifth-order RK method requires at least six stages. This naturally raises the question: Can the classical fourth-order explicit RK method be modified to achieve fifth-order accuracy without increasing the number of stages? The answer is yes. Incorporating radial basis function (RBF) interpolation into the method is one such approach. 

While traditional numerical methods for differential equations primarily rely on function evaluations, the availability of derivatives of solutions up to a certain order (whether known or estimated with sufficient accuracy) can significantly enhance the local accuracy of these methods. However, polynomial interpolation-based methods are inherently limited in their ability to incorporate such information in a flexible and adaptive manner. This limitation motivates the development of frameworks that can naturally accommodate derivative constraints. Such a framework can be obtained by implimentation of RBFs. An RBF $\phi$ is a function that is symmetric about an associated center $\mu_c$, i.e., $\phi_c(x)=\phi(\|x-\mu_c\|)$, where $\|\cdot\|$ denotes a vector norm \cite{Buhmann}. Common types of RBFs include the Gaussian RBF, multiquadratic (MQ) RBF, inverse-quadratic (IQ) RBF, and inverse-multiquadratic (IMQ) RBF, associated with a node $x_i$, are respectively defined as 
$\phi_i(x) =e^{-\varepsilon_i^2(x-x_i)^2}, \ \sqrt{1+\varepsilon_i^2(x-x_i)^2}, \ \left(1+\varepsilon_i^2(x-x_i)^2\right)^{-1} \ \text{and} \ \left(\sqrt{1+\varepsilon_i^2(x-x_i)^2}\right)^{-1},$
where $\varepsilon_i$ is the shape parameter. 
In RBF interpolation, the shape parameter is treated as a free variable that can adapt the local behavior of the solution. By varying this parameter, the accuracy of the interpolation can be significantly improved. Optimal values of the shape parameter are typically determined by formulating and solving an optimization problem, where the objective is to minimize the leading term of the truncation error.
Integrating RBFs into classical numerical schemes can reliably improve their accuracy while keeping them simple, which is especially beneficial when we have derivative information. 
The effectiveness of RBF techniques in enhancing the order and accuracy of finite-difference methods for first-order IVPs has been well demonstrated. In \cite{Gu20}, MQ-RBFs have been used to modify finite difference methods. Subsequently, the RBF-based Euler, midpoint, Adams-Bashforth, and Adams-Moultan methods have been developed using Gaussian and multiquadratic RBF interpolation in \cite{Gu21}. All of these RBF-based approaches outperform their classical counterparts in both accuracy and convergence order. Recently, in \cite{Gu25}, explicit RK methods have been generalized to explicit RBF RK methods by using the Gaussian RBF-based Euler method instead of the traditional Euler method, which raises the formal order by one. In the comparative study of the Gaussian-RBF-based and MQ-RBF-based Euler, midpoint, Adams-Bashforth, and Adams-Moultan methods, the authors in \cite{Gu21} observed that the superior performance of each RBF method depends on the considered problem; however, the adaptive MQ method demonstrates greater stability regions than the adaptive Gaussian method. Motivated by this, in this study, we employ the MQ-RBF Euler method to construct explicit MQ-RBF-based RK methods of order $s$ ($s=2,3,4$). This construction involves optimizing the shape parameters, which are $s-1$ in numbers when the RK method has $s$ stages. These parameters are locally optimized to reduce the local truncation error, thereby decreasing the global error and enhancing the local convergence order of the method. It is important to note that polynomial interpolation is a limiting case of MQ-RBF interpolation. Consequently, the MQ-RBF-based RK methods reduce to the classical RK methods when the shape parameters vanish. This ensures that the MQ-RBF-based RK methods maintain at least the same rate of convergence as the classical RK methods, even when the shape parameters deviate from their optimal values. 

It is well known that the accuracy of an RBF scheme improves as the shape parameter decreases, while instability increases exponentially with smaller shape parameter values \cite{Fornberg}. Optimal shape parameters refer to those values that strike a balance between accuracy and stability. While one may use trial-and-error or more systematic approaches to determine suitable values, there is no complete theoretical framework for selecting the “best” shape parameter, as its optimal value often depends on the problem under consideration \cite{Gu20}. Our approach to optimization is straightforward. Since truncation error analysis plays a crucial role in the numerical solution of IVPs, we use it to determine the optimal shape parameter. For an $s$-stage MQ-RBF-based RK method, the best shape parameter is the one that eliminates the leading truncation error term, i.e., the term of $\mathcal{O}(h^s)$ . Such a value is uniquely determined through this analysis. In this sense, our adaptive condition can be regarded as optimal. To derive the adaptive conditions for enhancing the order of convergence, we assume that the solution is sufficiently smooth. When this optimal shape parameter is employed in the proposed methods, tuned at each iteration, the overall accuracy of the method improves. As a result, the MQ-RBF-based RK methods become more efficient and accurate, with only a small increase in computational effort compared to the classical RK methods. We also plot the stability regions of MQ-RBF-based RK methods. To numerically show their consistency, accuracy, and convergence, we consider several specific IVPs, including stiff and non-stiff cases, and compare the results with their classical RK counterparts. Although the performance may be problem-specific, our results show that most $s$-stage MQ-RBF-based RK methods achieve better accuracy, a higher order of convergence, or both, compared to the classical $(s+1)$-stage RK methods. 

This article is organized as follows. In Section \ref{sec2}, we present the modification of $s$-stage classical explicit RK methods using adaptive MQ-RBF for $s=2,3,4$, and derive conditions on the shape parameters under which the local truncation error of the $s$-stage MQ-RBF-based RK method achieves order $s + 1$. Section \ref{sec3} provides a convergence analysis of the proposed MQ-RBF-based RK methods. In Section \ref{sec4}, we plot the stability regions of the MQ-RBF-based RK methods and compare them with those of the corresponding classical RK methods. Section \ref{sec5} presents numerical results that demonstrate the superiority of the MQ-RBF-based RK methods over their classical counterparts in terms of order and accuracy. Finally, Section \ref{sec6} concludes the article. 

\section{MQ-RBF-based RK methods}\label{sec2}
\noindent
Consider the initial value problem of the form
\begin{equation}\label{ivp}
u'(t)=f(t,u), \quad a\le t\le b,
\end{equation}
with initial condition
$$u(t_0)=u_0,$$
where we assume that $u(t) \in C^{\infty}[a, b]$ and $f(t, u)$ is a class of $C^\infty$ function. We divide the interval $[a, b]$ uniformly into subintervals $\{[t_i,t_{i+1}]\}_{i=1}^{N}$, where $t_i = a + ih, i = 0, 1, \ldots, N$, and $h = (b-a)/N$ is the size of the grid. 
We next derive the adaptive multiquadratic RBF based RK methods. Integrating \eqref{ivp} into the interval $[t_i,t_{i+1}]$, and applying the integral mean value theorem, we get
\begin{equation}\label{mvt}
u(t_{i+1})=u(t_i)+hf(t_i+\theta h, u(t_i+\theta h)), 
\end{equation}
where any value of $\theta \in [0,1]$ produces a numerical method.
When $\theta=0$, \eqref{mvt} yields
\begin{equation*}\label{euler}
u_{i+1}=u_i+hf(t_i, u_i),
\end{equation*}
which is the well known \textit{Euler} method.
When $\theta=1$, using \eqref{mvt} we obtain an implicit method as
\begin{equation}\label{th=1}
u_{i+1}=u_i+hf(t_{i+1}, u_{i+1}).
\end{equation}
Instead of using the classical \textit{Euler} method, as done in conventional RK methods to solve the problem explicitly, the idea is to approximate $u_{i+1}$ appearing on the right-hand side of \eqref{th=1} using the MQ-RBF \textit{Euler} method, i.e.,
\begin{equation}\label{rbfeu}
u_{i+1}=\left(1+\frac{\varepsilon_{2i}^2h^2}{2}\right)\left(u_i+hf(t_{i},u_i)\right).
\end{equation}
Applying \eqref{rbfeu} to \eqref{th=1}, we get
\begin{equation}
u_{i+1}=u_i+hf\left(t_{i+1}, \left(1+\frac{\varepsilon_{2i}^2h^2}{2}\right)u_i+h\left(1+\frac{\varepsilon_{2i}^2h^2}{2}\right)f(t_{i},u_i)\right).
\end{equation}
If we set 
\begin{align*}
K_1&=f(t_{i},u_i)\\
K_2&=f\left(t_{i+1}, \left(1+\frac{\varepsilon_{2i}^2h^2}{2}\right)u_i+h\left(1+\frac{\varepsilon_{2i}^2h^2}{2}\right)K_1\right)
\end{align*}
we get the method as 
\begin{equation}
u_{i+1}=u_i+hK_2.
\end{equation}
However, for $\theta=1/2$, \eqref{mvt} yields
\begin{equation}\label{th=1b2}
u_{i+1}=u_i+hf\left(t_{i}+\frac{h}{2}, u\left(t_{i}+\frac{h}{2}\right)\right).
\end{equation}
Note that $\left(t_{i}+\frac{h}{2}\right)$ is not a nodal point. Approximating $u\left(t_{i}+\frac{h}{2}\right)$ directly by the MQ-RBF \textit{Euler} method with spacing $h/2$, we get the method as
\begin{equation}
u_{i+1}=u_i+hf\left(t_{i}+\frac{h}{2}, \left(1+\frac{\varepsilon_{2i}^2h^2}{8}\right)u_i+\frac{h}{2}\left(1+\frac{\varepsilon_{2i}^2h^2}{8}\right)f(t_{i},u_i)\right).
\end{equation}
However, if we use the approximation
$$u'\left(t_i+\frac{h}{2}\right)\approx \frac{1}{2}\left(u'(t_i)+f(t_{i+1},u(t_j+h))\right),$$
with the MQ-RBF \textit{Euler} method, we get
$$f\left(t_{i}+\frac{h}{2}, u\left(t_{i}+\frac{h}{2}\right)\right)\approx \frac{1}{2}\left(f(t_i,u_i)+f\left(t_{i+1},\left(1+\frac{\varepsilon_{2i}^2h^2}{2}\right)\left(u_i+hf(t_{i},u_i)\right)\right)\right).$$
Thus, \eqref{th=1b2} can be approximated by
\begin{equation}
u_{i+1}=u_i+\frac{h}{2}\left(f(t_i,u_i)+f\left(t_{i+1},\left(1+\frac{\varepsilon_{2i}^2h^2}{2}\right)\left(u_i+hf(t_{i},u_i)\right)\right)\right),
\end{equation}
which can be written as 
\begin{equation}
u_{i+1}=u_i+\frac{h}{2}[K_1+K_2],
\end{equation}
where $K_1$ and $K_2$ is defined earlier.\\
When integrating \eqref{ivp} into the interval $[t_i,t_{i+1}]$, $f(t,u)$ can be considered as the slope of the solution curve that continuously varies over $[t_i,t_{i+1}]$. We may also interpret the methods in \eqref{th=1} and \eqref{th=1b2} as the cases when the slope of the solution curve in the interval $[t_i,t_{i+1}]$ is approximated by a single slope at the terminal point and midpoint, respectively. However, similar to RK methods,  the proposed method uses a weighted average of slopes over the given interval, rather than a single slope. Thus, we define the general MQ-RBF-based RK methods as
$$u_{i+1}=u_{i}+h[\text{weighted average of slopes over the given interval}].$$
Next, we derive explicit MQ-RBF-based RK methods with different stages. To determine the parameters involved, we expand $u_{i+1}$ and the function $f$ in powers of $h$, and optimize the value of the shape parameter $\varepsilon_i$, ensuring that the expansion matches the Taylor series of the exact solution up to a certain number of terms. This allows the proposed methods to achieve a higher order of convergence than conventional RK methods with the same number of stages.

\subsection{MQ-RBF-based RK2 method: Two-stage method of third-order accuracy}\label{MR1}
\noindent
Consider the following explicit MQ-RBF-based RK method with two slopes
\begin{align}\label{2me1}
\nonumber   K_1&=f(t_{i},\textbf{u}_i),\\
\nonumber   K_2&=f\left(t_{i}+c_2 h, \left(1+\frac{\varepsilon_{2i}^2(c_{2}h)^2}{2}\right)\left(\textbf{u}_i+ha_{21}K_1\right)\right),\\
\textbf{u}_{i+1}&=\textbf{u}_i+h[w_1K_1+w_2K_2],
\end{align}
where the parameters $c_2$, $a_{21}$, $w_1$ and $w_2$ selected to ensure that $\textbf{u}_{i+1}$ provides a close approximation to $\textbf{u}(t_{i+1})$.\\
The above two-stage method can be rewritten as the following iterative scheme:
\begin{equation}\label{2me2}
\begin{cases}
	\textbf{u}^{(1)}=\left(1+\frac{\varepsilon_{2i}^2(c_{2}h)^2}{2}\right)\textbf{u}_i+a_{21}h\left(1+\frac{\varepsilon_{2i}^2(c_{2}h)^2}{2}\right)f_i,\\
	\textbf{u}_{i+1}=\left[1-\frac{w_1}{a_{21}}\right]\textbf{u}_i+\frac{w_1}{a_{21}}\left(1+\frac{\varepsilon_{2i}^2(c_{2}h)^2}{2}\right)^{-1}\textbf{u}^{(1)}+w_2 h f(t_i+c_2h,\textbf{u}^{(1)}).
\end{cases}
\end{equation}
Let $u_i$ be the exact solution of the given initial value problem (IVP). Then, using the Taylor series expansion about $t_i$, the local truncation error can be obtained as
\begin{eqnarray}\label{e1}
\nonumber  \tau_i&=&\frac{\textbf{u}_{i+1}-\textbf{u}_{i}}{h}-(w_1K_1+w_2K_2)\\
\nonumber   &=&(1-w_1-w_2)f+h\left[\left(\frac{1}{2}-w_2c_2\right)f_t+\left(\frac{1}{2}-w_2a_{21}\right)f_uf\right]+h^2\bigg[\left(\frac{1}{6}-\frac{w_2c_2^2}{2}\right)f_{tt}\\
&&+\left(\frac{1}{3}-w_2c_2a_{21}\right)f_{tu}f+\left(\frac{1}{6}-\frac{w_2a_{21}^2}{2}\right)f^2f_{uu}+\left(\frac{f_t+ff_u}{6}-\frac{w_2c_{2}^2\varepsilon_{2i}^2}{2}u_i\right)f_u\bigg]+\mathcal{O}(h^3).
\end{eqnarray}
Consider 
\begin{equation}\label{p1}
a_{21}=c_2, \quad w_1=1-\frac{1}{2c_2} \quad \text{and} \quad w_2=\frac{1}{2c_2}, \quad c_2\ne 0.
\end{equation}
With these conditions on parameters and $\varepsilon_{2i}=0$, it is evident that the proposed scheme \eqref{2me2} reduces to the conventional RK2 method with second-order convergence.\\
Together with \eqref{p1}, if we further impose the condition
\begin{equation}\label{p2}
3w_2a_{21}^2=1,
\end{equation}
the parameters are then uniquely determined as
$$a_{21}=c_2=\frac{2}{3}, \quad w_1=\frac{1}{4} \quad \text{and} \quad w_2=\frac{3}{4},$$
causing the first three terms in the coefficient of $h^2$ to vanish. If we set $\varepsilon_{2i}^2=0$, then \eqref{2me2} reduces to \textit{Ralston}'s method, a second-order method with minimal local error represented by the remaining term. However, setting the remaining terms in the coefficient of $h^2$ to zero allows us to optimize $\varepsilon_{2i}$ as
\begin{equation}
\varepsilon_{2i}^2=\frac{f_t+ff_u}{u_i}=\frac{u_i''}{{u_i}},
\end{equation}
thereby elevating the scheme \eqref{2me2} to third-order accuracy.

\subsection{MQ-RBF-based RK3 method: Three-stage method of fourth-order accuracy}\label{MR2}
\noindent
We now use three evaluations of $f$ and define the MQ-RBF-based RK3 method as
\begin{equation}\label{3me1}
\textbf{u}_{i+1}=\textbf{u}_i+h[w_1K_1+w_2K_2+w_3K_3],
\end{equation}
with
\begin{align*}
K_1&=f(t_{i},\textbf{u}_i),\\
K_2&=f\left(t_{i}+c_2 h, \left(1+\frac{\varepsilon_{2i}^2(c_{2}h)^2}{2}\right)\left(\textbf{u}_i+ha_{21}K_1\right)\right),\\
K_3&=f\left(t_{i}+c_3 h, \left(1+\frac{\varepsilon_{3i}^2(c_{3}h)^2}{2}\right)\big(\textbf{u}_i+h(a_{31}K_1+a_{32}K_2)\big)\right).
\end{align*}
The alternative form of the above scheme is
\begin{equation}\label{3me2}
\begin{cases}
	\textbf{u}^{(1)}=\left(1+\frac{\varepsilon_{2i}^2(c_{2}h)^2}{2}\right)u_i+a_{21}h\left(1+\frac{\varepsilon_{2i}^2(c_{2}h)^2}{2}\right)f_i,\\
	\textbf{u}^{(2)}=\left(1+\frac{\varepsilon_{3i}^2(c_{3}h)^2}{2}\right)\left[\left(1-\frac{a_{31}}{a_{21}}\right)\textbf{u}_i+\frac{a_{31}}{a_{21}}\left(1+\frac{\varepsilon_{2i}^2(c_{2}h)^2}{2}\right)^{-1}\textbf{u}^{(1)}+a_{32}hf\left(t_i+c_2h,\textbf{u}^{(1)}\right)\right],\\
	\textbf{u}_{i+1}=\left[1-\frac{w_1}{a_{21}}-\frac{w_2}{a_{32}}\left(1-\frac{a_{31}}{a_{21}}\right)\right]\textbf{u}_i+\left[\frac{w_1}{a_{21}}-\frac{a_{31}w_2}{a_{32}a_{21}}\right]\left(1+\frac{\varepsilon_{2i}^2(c_{2}h)^2}{2}\right)^{-1}\textbf{u}^{(1)}\\
	\quad \quad \quad  +\frac{w_2}{a_{32}}\left(1+\frac{\varepsilon_{3i}^2(c_{3}h)^2}{2}\right)^{-1}\textbf{u}^{(2)}+w_3h f(t_i+c_3h,\textbf{u}^{(2)}).
\end{cases}
\end{equation}
By expanding $K_2$ and $K_3$ in the Taylor series about $t_i$ and substituting into \eqref{3me1}, and applying two conditions from the RK3 method given as
\begin{equation}\label{pv31}
a_{21}=c_2 \quad \text{and} \quad  a_{31}+a_{32}=c_3,
\end{equation}
the local truncation error can be simplified as
\begin{eqnarray}\label{e3}
\nonumber  \tau_i&=&\frac{\textbf{u}_{i+1}-\textbf{u}_{i}}{h}-(w_1K_1+w_2K_2+w_3K_3)\\
\nonumber   &=&(1-w_1-w_2-w_3)f+h\left[\left(\frac{1}{2}-w_2c_2-w_3c_3\right)(f_t+f_uf)\right]+h^2\bigg[\left(\frac{1}{6}-\frac{w_2c_2^2+w_3c_3^2}{2}\right)(f_{tt}+f^2f_{uu})\\
\nonumber  &&+\left(\frac{1}{3}-w_2c_2^2-w_3c_3^2\right)f_{tu}f+\left(\frac{1}{6}-\frac{a_{32}c_2w_3}{2}\right)(f_tf_u+f_u^2f)-\left(\frac{w_2c_{2}^2\varepsilon_{2i}^2+w_3c_{3}^2\varepsilon_{3i}^2}{2}\right)f_uu_i\bigg]\\
\nonumber&&+h^3\bigg[\left(\frac{1}{24}-\frac{w_2c_2^3+w_3c_3^3}{6}\right)(f_{ttt}+f_{uuu}f^3)+\left(\frac{1}{8}-a_{32}c_2c_3w_3\right)(f_tf_{tu}+f_tf_{uu}f)+\left(\frac{1}{24}-\frac{a_{32}c_2^2w_3}{2}\right)f_{tt}f_u\\
\nonumber &&\left(\frac{5}{24}-a_{32}c_2(c_2+c_3)w_3\right)f_{tu}f_uf+\left(\frac{1}{8}-\frac{w_2c_2^3+w_3c_3^3}{2}\right)(f_{ttu}f+f_{t uu}f^2)+\left(\frac{1}{6}-\frac{a_{32}c_2^2w_3}{2}-a_{32}c_2c_3w_3\right)f_{uu}f_uf^2\\
&&+\left(\frac{f_t+f_uf}{24}-\frac{a_{32}c_2^2w_3}{2}\varepsilon_{2i}^2u_i\right)f_u^2-\left(\frac{w_2c_2^3\varepsilon_{2i}^2+w_3c_3^3\varepsilon_{3i}^2}{2}\right)\left(u_i(f_{uu}f+f_{tu}) { +f_uf}\right) \bigg]+\mathcal{O}(h^4).
\end{eqnarray}
Clearly, applying the conditions
\begin{equation}\label{pv32}
w_1+w_2+w_3=1, \quad  2(w_2c_2+w_3c_3 )=1, \quad 3(w_2c_2^2+w_3c_3^2 )=1 \quad \text{and} \quad 6a_{32}c_2w_3=1,
\end{equation}
along with
\begin{equation}\label{pv33}
w_2c_2^2\varepsilon_{2i}^2+w_3c_3^2\varepsilon_{3i}^2=0,
\end{equation}
eliminates the constant term, $h$ term and $h^2$ term in the local truncation error \eqref{e3}, ensuring that the scheme is of order 3. Furthermore, if $\varepsilon_{2i}=\varepsilon_{3i}=0$, we recover the classical RK3 method. However, to achieve fourth-order accuracy, we need to optimize the shape parameters. For this, two additional conditions are required, which are specified in the following four distinct cases:
\begin{description}
\item[B1] $w_2c_2^3+w_3c_3^3=\frac{1}{4} \quad \text{and} \quad a_{32}c_2^2w_3=\frac{1}{12}$.\\
In this case the parameters are uniquely determined as
$$a_{21}=c_{2}=\frac{1}{2}, \quad a_{31}=-1, \quad a_{32}=2, \quad c_3=1, \quad w_1=\frac{1}{6}, \quad w_2=\frac{2}{3} \quad \text{and} \quad w_3=\frac{1}{6}.$$
In addition to the constant, $h$ and $h^2$ terms, these parameter values eliminate the coefficients of $f_{ttt}, \ f_{uuu}f^3, \ f_{tt}f_u$, $f_{ttu}f$ and $f_{t uu}f^2 $ in the $h^3$ term of \eqref{e3}. Now, by setting the remaining $h^3$ terms to zero, and using condition \eqref{pv32}, we obtain
$$\varepsilon_{2i}^2=\frac{{ (f_{uu}f-f_u^2+f_{tu})}u_i''}{{ (f_{uu}f-f_u^2+f_{tu})}u_i+{ f_uf}} \quad \text{and} \quad \varepsilon_{3i}^2=-\varepsilon_{2i}^2,$$
thereby making the proposed method fourth-order accurate. Note that this is only possible when the function $f$ satisfies condition $(f_{uu}f-f_u^2+f_{tu})u+f_uf\ne 0$ whenever $(f_{uu}f-f_u^2+f_{tu})(f_t+f_u f)\ne0$, otherwise the method exhibits only third-order accuracy.

\item[B2] $w_2c_2^3+w_3c_3^3=\frac{1}{4} \quad \text{and} \quad a_{32}c_2(c_2+c_3)w_3=\frac{5}{24}$.\\
For this case, we have two sets of solutions which are discussed below.
\begin{description}
	\item[(a)]  $a_{21}=c_{2}=\frac{5}{8}+\frac{\sqrt{33}}{24}, \quad a_{31}=-\frac{49}{256}+\frac{29\sqrt{33}}{768}, \quad a_{32}=\frac{209}{256}-\frac{61\sqrt{33}}{768}, \quad c_3=\frac{5}{8}-\frac{\sqrt{33}}{24}, \quad w_1=\frac{1}{8}, \quad w_2=\frac{7}{16}-\frac{3\sqrt{33}}{176} \quad \text{and} \quad w_3=\frac{7}{16}+\frac{3\sqrt{33}}{176}$. 
	These specific values eliminate the  $f_{ttt}, \ f_{uuu}f^3, \ f_{tu}f_uf, \ f_{ttu}f$ and $f_{tuu}f^2 $ terms from the coefficients of $h^3$ in \eqref{e3}.
	To eliminate the remaining  leading error term of $\mathcal{O}(h^3)$ we require 
	$$\varepsilon_{2i}^2=\dfrac{12f_u^2u_i''+(3+\sqrt{33})(f^2f_{uu}-f_{tt})f_u+2(3+\sqrt{33})(f_{uu}f+f_{tu})f_t}{\left[2(3+\sqrt{33})(f f_{uu}+f_{tu})+(15+\sqrt{33})f_u^2\right]u_i{ +2(3+\sqrt{33})f_uf}} \ \ \text{and} \ \ \varepsilon_{3i}^2=\frac{-7-\sqrt{33}}{4}\varepsilon_{2i}^2.$$
	
	\item[(b)]  $a_{21}=c_{2}=\frac{5}{8}-\frac{\sqrt{33}}{24}, \quad a_{31}=-\frac{49}{256}-\frac{29\sqrt{33}}{768}, \quad a_{32}=\frac{209}{256}+\frac{61\sqrt{33}}{768}, \quad c_3=\frac{5}{8}+\frac{\sqrt{33}}{24}, \quad w_1=\frac{1}{8}, \quad w_2=\frac{7}{16}+\frac{3\sqrt{33}}{176} \quad \text{and} \quad w_3=\frac{7}{16}-\frac{3\sqrt{33}}{176}$.
	Furthermore, by eliminating the leading error term of $\mathcal{O}(h^3)$ in \eqref{e3}, we get 
	$$\varepsilon_{2i}^2=\dfrac{12f_u^2u_i''+(3-\sqrt{33})(f_{uu}f^2-f_{tt})f_u+2(3-\sqrt{33})(f_{uu}f+f_{tu})f_t}{\left[2(3-\sqrt{33})(f f_{uu}+f_{tu})+(15-\sqrt{33})f_u^2\right]u_i{ +2(3-\sqrt{33})f_uf}} \ \ \text{and} \ \ \varepsilon_{3i}^2=\frac{-7+\sqrt{33}}{4}\varepsilon_{2i}^2.$$
	
\end{description}
\item[B3] $w_2c_2^3+w_3c_3^3=\frac{1}{4} \quad \text{and} \quad \frac{a_{32}c_2^2w_3}{2}+a_{32}c_2c_3w_3=\frac{1}{6}$.\\
These conditions eliminate the terms $f_{ttt}, \ f_{uuu}f^3, \  f_{ttu}f, \ f_{tuu}f^2$  and $f_{uu}f_uf^2$ from the coefficients of $h^3$ in \eqref{e3}, and provide two sets of parameter values.
\begin{description}
	\item[(a)] $a_{21}=c_{2}= 1, \quad a_{31}=\frac{1}{4}, \quad a_{32}=\frac{1}{4}, \quad c_3= \frac{1}{2}, \quad w_1=\frac{1}{6}, \quad w_2=\frac{1}{6}  \quad \text{and} \quad w_3= \frac{2}{3}$. 
	Further optimizing the shape function as
	$$\varepsilon_{2i}^2=\dfrac{f_u^2u_i''-(f_{tu}f+f_{tt})f_u+(f_{uu}f+f_{tu})f_t}{\left[f_{uu}f+f_{tu}+2f_u^2\right]u_i{ +f_uf}} \ \ \text{and} \ \ \varepsilon_{3i}^2=-\varepsilon_{2i}^2,$$
	We obtain the MQ-RBF-based RK3 method of fourth-order accuracy.
	\item[(b)] $a_{21}=c_{2}= \frac{1}{3}, \quad a_{31}=-\frac{5}{12}, \quad a_{32}=\frac{5}{4}, \quad c_3= \frac{5}{6}, \quad w_1=\frac{1}{10}, \quad w_2=\frac{1}{2}  \quad \text{and} \quad w_3= \frac{2}{5}$. We further optimize the shape parameters as 
	$$\varepsilon_{2i}^2=\dfrac{3f_u^2u_i''+(f_{tu}f+f_{tt})f_u-(f_{uu}f+f_{tu})f_t}{\left[-f_{uu}f-f_{tu}+2f_u^2\right]u_i{ -f_uf}} \ \ \text{and} \ \ \varepsilon_{3i}^2=-\frac{1}{5}\varepsilon_{2i}^2,$$
	to obtain the fourth-order MQ-RBF-based RK3 method.
\end{description}
\item[B4] $a_{32}c_2c_3w_3=\frac{1}{8} \quad \text{and} \quad a_{32}c_2^2w_3=\frac{1}{12}$.\\
These conditions provide unique values for parameters used as
$a_{21}=c_{2}= \frac{1}{2}, \quad a_{31}=0, \quad a_{32}=\frac{3}{4}, \quad c_3= \frac{3}{4}, \quad w_1=\frac{2}{9}, \quad w_2=\frac{1}{3}  \quad \text{and} \quad w_3= \frac{4}{9}$,
which annihilate the terms $f_{t}f_{tu}$, $f_{t}f_{uu}f$, $f_{tt}f_u$, $f_{tu}f_u f$ and $f_{uu}f_uf^2$ from the coefficient of $h^3$ in \eqref{e3}. We then equate the remaining terms in the coefficient of $h^3$ in \eqref{e3} to zero in order to obtain the shape functions as
$$\varepsilon_{2i}^2=\dfrac{12f_u^2u_i''+f_{ttt}+f_{uuu}f^3+3(f_{ttu}+f_{tuu}f)f}{3\left(-f_{uu}f-f_{tu}+4f_u^2\right)u_i{ -3f_uf}} \ \ \text{and} \ \ \varepsilon_{3i}^2=-\frac{1}{3}\varepsilon_{2i}^2.$$
Evidently, scheme \eqref{3me2} along with given parameters is of fourth-order accuracy.
\end{description}
\subsection{MQ-RBF-based RK4 method: Four-stage method of fifth-order accuracy} \label{MR3}
\noindent 
It is evident from the scheme presented in the preceding subsection that including a shape parameter at each stage can analytically enhance the order of the scheme by one. We apply the same technique to construct the four-stage MQ-RBF-based RK method as
\begin{equation}\label{4me1}
\textbf{u}_{i+1}=\textbf{u}_i+h[w_1K_1+w_2K_2+w_3K_3+w_4K_4],
\end{equation}
with
\begin{align*}
K_1&=f(t_{i},\textbf{u}_i),\\
K_2&=f\left(t_{i}+c_2 h, \left(1+\frac{\varepsilon_{2i}^2(c_{2}h)^2}{2}\right)\left(\textbf{u}_i+ha_{21}K_1\right)\right),\\
K_3&=f\left(t_{i}+c_3 h, \left(1+\frac{\varepsilon_{3i}^2(c_{3}h)^2}{2}\right)\big(\textbf{u}_i+h(a_{31}K_1+a_{32}K_2)\big)\right)\\
K_4&=f\left(t_{i}+c_4 h, \left(1+\frac{\varepsilon_{4i}^2(c_{4}h)^2}{2}\right)\big(\textbf{u}_i+h(a_{41}K_1+a_{42}K_2+a_{43}K_3)\big)\right),
\end{align*}
which can alternatively be written as
\begin{equation}\label{4me2}
\begin{cases}
	\textbf{u}^{(1)}=\left(1+\frac{\varepsilon_{2i}^2(c_{2}h)^2}{2}\right)\textbf{u}_i+a_{21}h\left(1+\frac{\varepsilon_{2i}^2(c_{2}h)^2}{2}\right)f_i,\\
	\textbf{u}^{(2)}=\left(1+\frac{\varepsilon_{3i}^2(c_{3}h)^2}{2}\right)\left[\left(1-\frac{a_{31}}{a_{21}}\right)\textbf{u}_i+\frac{a_{31}}{a_{21}}\left(1+\frac{\varepsilon_{2i}^2(c_{2}h)^2}{2}\right)^{-1}\textbf{u}^{(1)}+a_{32}hf\left(t_i+c_2h,\textbf{u}^{(1)}\right)\right],\\
	\textbf{u}^{(3)}=\left(1+\frac{\varepsilon_{4i}^2(c_{4}h)^2}{2}\right)\Big[\left(1-\frac{a_{41}}{a_{21}}-\frac{a_{42}}{a_{32}}\left(1-\frac{a_{31}}{a_{21}}\right)\right)\textbf{u}_i+\left(\frac{a_{41}}{a_{21}}-\frac{a_{31}a_{42}}{a_{32}a_{21}}\right)\left(1+\frac{\varepsilon_{2i}^2(c_{2}h)^2}{2}\right)^{-1}\textbf{u}^{(1)}\\
	\quad \quad \quad+\frac{a_{42}}{a_{32}} \left(1+\frac{\varepsilon_{3i}^2(c_{3}h)^2}{2}\right)^{-1}\textbf{u}^{(2)} +a_{43}hf\left(t_i+c_3h,\textbf{u}^{(2)}\right)\Big],\\
	\textbf{u}_{i+1}=\left[1-\frac{w_1}{a_{21}}-\frac{w_2}{a_{32}}\left(1-\frac{a_{31}}{a_{21}}\right)-\frac{w_3}{a_{43}}\left(1-\frac{a_{41}}{a_{21}}-\frac{a_{42}}{a_{32}}\left(1-\frac{a_{31}}{a_{21}}\right)\right)\right]\textbf{u}_i+\left[\frac{w_1}{a_{21}}-\frac{a_{31}w_2}{a_{32}a_{21}}-\frac{w_3}{a_{43}}\left(\frac{a_{41}}{a_{21}}-\frac{a_{31}a_{42}}{a_{32}a_{21}}\right)\right]\\
	\quad \quad \quad \left(1+\frac{\varepsilon_{2i}^2(c_{2}h)^2}{2}\right)^{-1}\textbf{u}^{(1)} +\left[\frac{w_2}{a_{32}}-\frac{w_3a_{42}}{a_{43}a_{32}} \right]\left(1+\frac{\varepsilon_{3i}^2(c_{3}h)^2}{2}\right)^{-1}\textbf{u}^{(2)}+\frac{w_3}{a_{43}}\left(1+\frac{\varepsilon_{4i}^2(c_{4}h)^2}{2}\right)^{-1}\textbf{u}^{(3)}\\
	\quad \quad \quad +w_4h f(t_i+c_4h,\textbf{u}^{(3)}).
\end{cases}
\end{equation}
Under the conditions
\begin{equation}\label{pv41}
a_{21}=c_2, \quad a_{31}+a_{32}=c_3, \quad a_{41}+a_{42}+a_{43}=c_4,
\end{equation}
the Taylor series expansion of \eqref{4me1} yields the local truncation error as
\begin{eqnarray}\label{e4}
\nonumber  \tau_i&=&\frac{\textbf{u}_{i+1}-\textbf{u}_{i}}{h}-(w_1K_1+w_2K_2+w_3K_3+w_4K_4)\\
\nonumber   &=&(1-w_1-w_2-w_3-w_4)f+h\left[\left(\frac{1}{2}-w_2c_2-w_3c_3-w_4c_4\right)(f_t+f_uf)\right]+h^2\bigg[\left(\frac{1}{6}-\frac{w_2c_2^2+w_3c_3^2+w_4c_4^2}{2}\right)\\ 
\nonumber  &&(f_{tt}+f^2f_{uu}+2f_{tu}f)
+\left(\frac{1}{6}-a_{32}c_2w_3+w_4(c_2a_{42}+c_3a_{43})\right)(f_tf_u+f_u^2f) -\left(\frac{w_2c_{2}^2\varepsilon_{2i}^2+w_3c_{3}^2\varepsilon_{3i}^2+w_4c_{4}^2\varepsilon_{4i}^2}{2}\right)f_uu_i\bigg]\\
\nonumber&&+h^3\bigg[\left(\frac{1}{24}-\frac{w_2c_2^3+w_3c_3^3+w_4c_4^3}{6}\right)(f_{ttt}+f_{uuu}f^3+3f_{ttu}f+3f_{t uu}f^2)\\
\nonumber &&+\left(\frac{1}{8}-a_{32}c_2c_3w_3-(a_{42}c_2+a_{43}c_3)c_4w_4\right)(f_tf_{tu}+f_tf_{uu}f)+\left(\frac{1}{24}-\frac{a_{32}c_2^2w_3+(a_{42}c_2^2+a_{43}c_3^2)w_4}{2}\right)f_{tt}f_u\\
\nonumber &&+\left(\frac{5}{24}-a_{32}c_2(c_2+c_3)w_3-a_{42}c_2(c_2+c_4)w_4-a_{43}c_3(c_3+c_4)w_4\right)f_{tu}f_uf\\
\nonumber&&+\left(\frac{1}{6}-a_{32}c_2\left(\frac{c_2}{2}+c_3\right)w_3-a_{42}c_2\left(\frac{c_2}{2}+c_4\right)w_4-a_{43}c_3\left(\frac{c_3}{2}+c_4\right)w_4\right)f_{uu}f_uf^2\\
\nonumber&&+\left(\left(\frac{1}{24}-a_{32}a_{43}c_2w_4\right)(f_t+ff_u)-\frac{(a_{32}w_3+a_{42}w_4)c_2^2\varepsilon_{2i}^2+a_{43}c_3^2w_4\varepsilon_{3i}^2}{2}u_i\right)f_u^2 \\
\nonumber&&-\left(\frac{w_2c_2^3\varepsilon_{2i}^2+w_3c_3^3\varepsilon_{3i}^2+w_4c_4^3\varepsilon_{4i}^2}{2}\right)\left(u_i(f_{uu}f+f_{tu}) { +f_uf}\right)\bigg]+h^4\bigg[\left(\frac{1}{120}-\frac{w_2c_2^4+w_3c_3^4+w_4c_4^4}{24}\right)(f_{tttt}+f_{uuuu}f^4\\
\nonumber&&+4f_{tttu}f+4f_{tuuu}f^3 +6f_{ttuu}f^2)+\left(\frac{1}{20}-\frac{a_{32}c_2c_3^2w_3+(a_{42}c_2+a_{43}c_3)c_4^2w_4}{2}\right)(f_{ttu}f_t+2f_{tuu}f_tf+f_{uuu}f_tf^2)\\
\nonumber&& +\left(\frac{1}{30}-\frac{a_{32}c_2^2c_3w_3+(a_{42}c_2^2+a_{43}c_3^2)c_4w_4}{2}\right)(f_{tt}f_{tu}f+2f_{tu}^2f+f_{tt}f_{uu}f+f_{uu}^2f^3+3f_{tu}f_{uu}f^2)\\
\nonumber&&+\left(\frac{1}{120}-\frac{a_{32}c_2^3w_3+(a_{42}c_2^3+a_{43}c_3^3)w_4}{6}\right)f_{ttt}f_u+\left(\frac{7}{120}-a_{43}a_{32}c_2(c_3+c_4)w_4\right)f_tf_{tu}f_u+\left(\frac{1}{120}-\frac{a_{32}a_{43}c_2^2w_4}{2}\right)f_{tt}f_u^2\\
\nonumber&&+\left(\frac{1}{40}-\frac{a_{32}^2c_2^2w_3+(a_{42}c_2+a_{43}c_3)^2w_4}{2}\right)f_{uu}f_{t}^2+\left(\frac{3}{40}-\frac{a_{32}c_2(c_2^2+c_3^2)w_3+a_{42}c_2(c_2^2+c_4^2)w_4+a_{43}c_3(c_3^2+c_4^2)w_4}{2}\right)\\
\nonumber&&f_{ttu}f_{u}f+\left(\frac{1}{8}-\frac{a_{32}c_2(c_2^2+2c_3^2)w_3+a_{42}c_2(c_2^2+2c_4^2)w_4+a_{43}c_3(c_3^2+2c_4^2)w_4}{2}\right)f_{tuu}f_uf^2 \\
\nonumber&&+\left(\frac{3}{40}-a_{32}a_{43}c_2(c_2+c_3+c_4)w_4\right)f_{tu}f_{u}^2f+\bigg(\frac{13}{120}-a_{32}^2c_2^2w_3-a_{32}a_{43}c_2(c_3+c_4)w_4-a_{42}a_{43}c_2c_3w_4-a_{42}^2c_2^2w_4\\
\nonumber&& -a_{43}^2c_3^2w_4\bigg)f_{uu}f_tf_uf+\bigg(\frac{11}{120}-\frac{a_{32}^2c_2^2w_3}{2}-\frac{2a_{32}a_{43}c_2(c_3+c_4)+2a_{42}a_{43}c_2c_3+a_{42}^2c_2^2 +a_{43}^2c_3^2+a_{32}a_{43}c_2^2}{2}w_4\bigg)f_{uu}f_u^2f^2\\
\nonumber&&+\left(\frac{7}{120}-\frac{a_{32}c_2(c_2^2+3c_3^2)w_3+(a_{42}c_2(c_2^2+3c_4^2)+a_{43}c_3(c_3^2+3c_4^2))w_4}{6}\right)f_{uuu}f_uf^3-\frac{(w_2c_2^4\varepsilon_{2i}^2+w_3c_3^4\varepsilon_{3i}^2+w_4c_4^4\varepsilon_{4i}^2)}{4}\\
\nonumber&& \left(\left(f_{ttu}+2f_{tuu}+f_{uuu}\right)u_{i} { +2(f_{uu}f^2+f_{tu}f)}\right){ -\frac{a_{32}w_3c_2^3\varepsilon_{2i}^2+w_4a_{42}c_2^3\varepsilon_{2i}^2+w_4a_{43}c_3^3\varepsilon_{3i}^2}{2}f_u^2f}\\
\nonumber &&-\frac{w_3a_{32}(c_2+c_3)c_2^2\varepsilon_{2i}^2+w_4a_{42}(c_2+c_4)c_2^2\varepsilon_{2i}^2+w_4a_{43}(c_3+c_4)c_3^2\varepsilon_{3i}^2}{2}\left(f_{tu}f_u+f_{uu}f_uf\right)u_i+\frac{a_{32}a_{43}c_2^2w_4\varepsilon_{2i}^2}{2}f_u^3u_i \\
\nonumber && -\frac{a_{32}c_2w_3c_3^2\varepsilon_{3i}^2+(a_{42}c_2+a_{43}c_3)w_4c_4^2\varepsilon_{4i}^2}{2}\left(\left(f_tf_{uu}+f_uf_{uu}f\right)u_i+{ f_u^2f+f_tf_u}\right)-\frac{(w_2c_2^4\varepsilon_{2i}^4+w_3c_3^4\varepsilon_{3i}^4+w_4c_4^4\varepsilon_{4i}^4)}{8}{ f_{uu}u_i^2}\\
&&+\frac{f_u^3u_i''}{120}\bigg]+\mathcal{O}(h^5).
\end{eqnarray}
If we impose the conditions 
\begin{eqnarray}\label{pv42}
\nonumber  & w_1+w_2+w_3+w_4=1, \quad   2(w_2c_2+w_3c_3+w_4c_4)=1, \quad 3(w_2c_2^2+w_3c_3^2+w_4c_4^2)=1, \quad 4(w_2c_2^3+w_3c_3^3+w_3c_4^3)=1\\
\nonumber &6(a_{32}w_3c_2+a_{42}w_4c_2+a_{43}w_4c_3)=1, \quad 8(a_{32}w_3c_2c_3+a_{42}w_4c_2c_4+a_{43}w_4c_3c_4)=1, \quad 12(a_{32}w_3c_2^2+a_{42}w_4c_2^2+a_{43}w_4c_3^2)=1,\\
&24a_{32}a_{43}w_4c_2=1,
\end{eqnarray}
along with 
\begin{eqnarray}\label{pc43}
\nonumber &w_2c_2^2\varepsilon_{2i}^2+w_3c_3^2\varepsilon_{3i}^2+w_4c_4^2\varepsilon_{4i}^2=0, \quad w_2c_2^3\varepsilon_{2i}^2+w_3c_3^3\varepsilon_{3i}^2+w_4c_4^3\varepsilon_{4i}^2=0, \quad a_{32}w_3c_2^2\varepsilon_{2i}^2+a_{42}w_4c_2^2\varepsilon_{2i}^2+a_{43}w_4c_3^2\varepsilon_{3i}^2=0, \\
\end{eqnarray}
then all the constant, $h$, $h^2$, and $h^3$ terms in $\tau_n$ vanish, indicating that our scheme, namely the MQ-RBF-based RK4 method, is fourth-order accurate. However, to achieve fifth-order accuracy, the $h^4$ terms must also vanish. To this end, we propose two possible values for the shape parameters.
\begin{description}
\item[C1] Consider $a_{32}c_2c_3^2w_3+(a_{42}c_2+a_{43}c_3)c_4^2w_4=\frac{1}{10}$ and $a_{32}c_2^2c_3w_3+(a_{42}c_2^2+a_{43}c_3^2)c_4w_4=\frac{1}{15}$, so that we obtain
\begin{eqnarray*}
	&a_{21}=c_2=\frac{2}{5}, \quad a_{31}=-\frac{3}{20}, \quad a_{32}=\frac{3}{4}, \quad a_{41}=\frac{19}{44},\quad a_{42}=-\frac{15}{44}, \quad a_{43}=\frac{10}{11}, \quad c_{3}=\frac{3}{5}, \quad c_{4}=1,\\
	&w_{1}=w_4=\frac{11}{72}, \quad w_{2}=w_3=\frac{25}{72}.
\end{eqnarray*}
Now, using \eqref{pc43} and equating the leading error term in $\mathcal{O}(h^4)$ to zero, we get
$$ \varepsilon_{3i}^2=-\frac{2}{3}\varepsilon_{2i}^2, \quad \varepsilon_{4i}^2=\frac{2}{11}\varepsilon_{2i}^2,$$
where $\varepsilon_{2i}^2$ is a real valued solution of the equation $\alpha x^2+\beta x+\gamma=0 $, with 
\begin{eqnarray*}
	\alpha&=&168 f_{uu}u_i^2\\
	\beta&=& (  66f_{ttu}+66f_{uuu} f^2+ 132f_{tuu}f- 462f_{tu}f_u- 270f_{uu}f_t-732f_{uu}f_uf+ 330f_u^3)u_i\\
	&& { + 132f^2f_{uu}-402f_u^2f+132f_{tu}f-270f_tf_u}\\
	\gamma&=&11(f_{tttt}+f_{uuuu}f^4+4 f_{tttu}f+4f_{tuuu}f^3+6f_{ttuu}f^2)- 44 (f_{ttt}f_u+3f_{ttu}f_u f+3f_u f_{tuu} f^2+f_u  f_{uuu}f^3)\\
	&&+ 330 f_t f_{tu} f_u +330 f_{tu} f_u^2f+ 135 f_t^2  f_{uu}+600 f_t  f_uf_{uu}f +465f_u^2  f_{uu}f^2-330 f_u^3u_i''. 
\end{eqnarray*}
\item[C2] Consider $a_{32}c_2c_3^2w_3+(a_{42}c_2+a_{43}c_3)c_4^2w_4=\frac{1}{10}$ and $a_{32}^2c_2^2w_3+(a_{42}c_2+a_{43}c_3)^2w_4=\frac{1}{20}$, so that we obtain
\begin{eqnarray*}
	&a_{21}=c_2=\frac{1}{4}, \quad a_{31}=-\frac{6}{25}, \quad a_{32}=\frac{21}{25}, \quad a_{41}=\frac{6}{5},\quad a_{42}=-\frac{57}{35}, \quad a_{43}=\frac{10}{7}, \quad c_{3}=\frac{3}{5}, \quad c_{4}=1,\\
	&w_{1}=\frac{1}{9}, \quad w_{2}=\frac{16}{63}, \quad w_{3}=\frac{125}{252}, \quad w_{4}=\frac{5}{36}.
\end{eqnarray*}
Now, using \eqref{pc43} and equating the leading error term in $\mathcal{O}(h^4)$ to zero, we get
$$\varepsilon_{3i}^2=-\frac{\varepsilon_{2i}^2}{6}, \ \ \text{and}\ \ \varepsilon_{4i}^2=\frac{\varepsilon_{2i}^2}{10},$$
where $\varepsilon_{2i}^2$ is a real valued solution of the equation $\alpha x^2+\beta x+\gamma=0 $, 
with 
\begin{eqnarray*}
	\alpha&=&3f_{uu}u_i^2\\
	\beta&=&6 (f_{ttu}+f^2f_{uuu}+2f_{tuu}f-7f_{tu}f_u-7f_{uu}f_uf+5f_u^3)u_i {  +12(f_{uu}f^2-f_u^2f+f_{tu}f)}\\
	\gamma&=&f_{tttt}+ 4 f_{tttu}f+ 6 f_{ttuu}f^2+ 4 f_{tuuu}f^3+f_{uuuu}f^4 - 4(f_{ttt}+ 3 f_{ttu} f+3 f_{tuu} f^2+ f_{uuu}f^3)f_u\\
	&&+ 18f_{tt} f_{tu} +18f_{tt} f_{uu}f + 18 f_{uu}^2f^3+36 f_{tu}^2f+ 54f_{tu} f_{uu}f^2 + 30f_{uu}f_u^2 f^2 + 48 f_{uu} f_t f_u f + 48 f_{tu}f_t f_u   \\
	&& + 12 f_{tu} f_u^2f - 18f_{tt} f_u^2 - 48f_u^3u_i''
\end{eqnarray*}
\end{description}
It is important to note that in both cases described above, the resulting quadratic equations may not yield real roots. In such situations, eliminating the $h^4$ term in \eqref{e4} is not always feasible. Nevertheless, the MQ-RBF-based RK4 method \eqref{3me1} still achieves fourth-order accuracy.
\section{Convergence of MQ-RBF-based RK methods}\label{sec3}
\noindent
To analyze the convergence of the MQ-RBF-based RK methods proposed in Section~\ref{sec2}, we study the behavior of the difference between the exact solution of IVP \eqref{ivp} and the corresponding numerical solution. To guarantee the unique existence of an exact solution to \eqref{ivp}, we assume that $f(t, u)$ is Lipschitz continuous with a Lipschitz constant $L$ throughout this section. Suppose that the interval $[a,b]$ is partitioned into $N$ equal subintervals, resulting in a uniform grid of $N+1$ nodes, i.e., 
$$a=t_0<t_1<\cdots<t_N=b,$$
where for each $n=0,1,\ldots,N$, the grid nodes are defined as $t_n=a+nh$, with $h=(b-a)/N$. 
Let $u_i = u(t_i)$ denote the exact solution at time $t_i$, and let $\textbf{u}_i$ denote its numerical approximation obtained using the MQ-RBF-based RK methods. Then we have the following results:
\begin{theorem}\label{th1}
Consider that $\varepsilon_{2i}^2$ is bounded for all $i=0,1,\ldots,N-1$. Then the MQ-RBF-based RK2 method \eqref{2me1}, satisfying \eqref{p1}, converges.
\end{theorem}
\textit{Proof.} Let $E_{i+1}:=|u_{i+1}-\textbf{u}_{i+1}|$. For a fixed $t=t_i$, \eqref{2me1} yields
\begin{eqnarray*}
K_1(v)&=&f(t_{i},v),\\
K_2(v)&=&f\left(t_{i}+c_2 h, \left(1+\frac{\varepsilon_{2i}^2(c_{2}h)^2}{2}\right)\left(v+ha_{21}K_1(v)\right)\right).
\end{eqnarray*}
It is clear that the exact solution of \eqref{ivp} satisfies
\begin{equation}\label{2con1}
u_{i+1}=u_i+h(w_1K_1(u_i)+w_2K_2(u_i))+h\tau_i,
\end{equation}
where $\tau_i$ is the local truncation error defined in \eqref{e1}. Moreover, due to the Lipschitz continuity of $f(t,u)$ with respect to $u$, it follows that
$$|K_1(u_i)-K_1(\textbf{u}_i)|=|f(t_i,u_i)-f(t_i,\textbf{u}_i)|\le L|u_i-\textbf{u}_i|,$$ and
\begin{eqnarray*}
\Big|K_2(u_i)-K_2(\textbf{u}_i)|&=&|f\left(t_{i}+c_2 h, \left(1+\frac{\varepsilon_{2i}^2(c_{2}h)^2}{2}\right)\left(u_i+ha_{21}K_1\right)\right)-f\left(t_{i}+c_2 h, \left(1+\frac{\varepsilon_{2i}^2(c_{2}h)^2}{2}\right)\left(\textbf{u}_i+ha_{21}K_1\right)\right)\Big|\\
&\le& L\Big|\left(1+\frac{\varepsilon_{2i}^2(c_{2}h)^2}{2}\right)\left(u_i+ha_{21}K_1(u_i)\right)-\left(1+\frac{\varepsilon_{2i}^2(c_{2}h)^2}{2}\right)\left(\textbf{u}_i+ha_{21}K_1(\textbf{u}_i)\right)\Big|\\
&\le& L\left(1+\frac{\varepsilon_{2i}^2(c_{2}h)^2}{2}\right)\Big(|u_i-\textbf{u}_i|+ha_{21}|K_1(u_i)-K_1(\textbf{u}_i)|\Big)\\
&\le& L(1+ha_{21}L)\left(1+\frac{\varepsilon_{2i}^2(c_{2}h)^2}{2}\right)|u_i-\textbf{u}_i|.
\end{eqnarray*}
Now, subtracting \eqref{2con1} from \eqref{2me1}, we obtain
\begin{eqnarray}\label{ec1}
\nonumber   E_{i+1}&=&|u_i-\textbf{u}_i+hw_1(K_1(u_i)-K_1(\textbf{u}_i))+hw_2(K_2(u_i)-K_2(\textbf{u}_i))+h\tau_i|\\
&\le&\varphi_iE_i+h|\tau_i|,
\end{eqnarray}
where $\varphi_i=1+(w_1+w_2)hL+w_2a_{21}h^2L^2+\left(w_2hL+w_2a_{21}h^2L^2\right)\frac{c_2^2h^2\varepsilon_{2i}^2}{2}$, and $E_i=|u_i-\textbf{u}_i|$. Let 
$$\mathcal{C}_2=\displaystyle \max_{i=0,\ldots,N-1}\left(\sup_{0< h\le b-a}\frac{\left(2w_2L+hL^2\right)c_2^2h^2\varepsilon_{2i}^2}{4e^{hL}}\right).$$
Since $\varepsilon_i$ is bounded for all $i=0,1, \ldots, N-1$, we have $\mathcal{C}_2<\infty$. Thus, we have
\begin{eqnarray*}
\varphi_i&=&1+hL+\frac{1}{2}h^2L^2+\left(w_2hL+\frac{1}{2}h^2L^2\right)\frac{c_2^2h^2\varepsilon_{2i}^2}{2}\\
&\le&\left(1+\mathcal{C}_2 \right)e^{hL}\\
&\le&e^{(L+\mathcal{C}_2)h}.
\end{eqnarray*}
and consequently, by induction, \eqref{ec1} yields
\begin{equation}
E_i\le \left(e^{(L+\mathcal{C}_2)h}\right)^{i-1}\varphi_kE_0+h\sum_{k=1}^{i-1}\left(e^{(L+\mathcal{C}_2)h}\right)^{i-k}|\tau_{k-1}|+h|\tau_{i-1}|.
\end{equation}
Assuming $E_0=0$, we obtain
\begin{equation*}
E_i\le ihe^{(L+\mathcal{C}_2)ih}\|\tau\|_{\infty}\le (b-a)e^{(L+\mathcal{C}_2)(b-a)}\|\tau\|_{\infty},
\end{equation*}
for all $i=0,1,\ldots N-1$, where $\|\tau\|_{\infty}=\displaystyle \max_{0,\ldots,N-1}|\tau_i|$.
Therefore, for $h=(b-a)/N$, we have
$$\lim_{h\to 0}E_i=0.$$
This implies that $\textbf{u}_N\to u_N=u(b)$ as $h$ tends to zero. Hence, the MQ-RBF-based RK2 method \eqref{2me1} converges. 
\begin{theorem}
Consider that $\varepsilon_{2i}^2$ and $\varepsilon_{3i}^2$ are bounded for all $i=0,1,\ldots,N-1$. Then the MQ-RBF-based RK3 method \eqref{3me1}, satisfying \eqref{pv31} and \eqref{pv32}, is convergent.
\end{theorem}
\textit{Proof.} By the definition of local truncation error, from \eqref{e3} we can write
\begin{equation}\label{3con1}
u_{i+1}=u_i+h(w_1K_1(u_i)+w_2K_2(u_i)+w_3K_3(u_i))+h\tau_i.
\end{equation}
Let us fix $t=t_i$ and define
\begin{align*}
K_1(v)&=f(t_{i},v),\\
K_2(v)&=f\left(t_{i}+c_2 h, \left(1+\frac{\varepsilon_{2i}^2(c_{2}h)^2}{2}\right)\left(v+ha_{21}K_1(v)\right)\right),\\
K_3(v)&=f\left(t_{i}+c_3 h, \left(1+\frac{\varepsilon_{3i}^2(c_{3}h)^2}{2}\right)\big(v+h(a_{31}K_1(v)+a_{32}K_2(v))\big)\right),
\end{align*}
accordingly. Then, invoking the Lipschitz continuity of $f(t,u)$ with respect to $u$, we can estimate
\begin{eqnarray*}
|K_1(u_i)-K_1(\textbf{u}_i)|&\le& L|u_i-\textbf{u}_i|,\\
|K_2(u_i)-K_2(\textbf{u}_i)| &\le& L(1+ha_{21}L)\left(1+\frac{\varepsilon_{2i}^2(c_{2}h)^2}{2}\right)|u_i-\textbf{u}_i|,
\end{eqnarray*}
and
\begin{eqnarray*}
|K_3(u_i)-K_3(\textbf{u}_i)|&\le& L\left(1+\frac{\varepsilon_{3i}^2(c_{3}h)^2}{2}\right)\left(1+h(a_{31}+a_{32})L+a_{21}a_{32}h^2L^2+(1+hLa_{21})a_{32}hL\frac{\varepsilon_{2i}^2(c_{2}h)^2}{2}\right)|u_i-\textbf{u}_i|.
\end{eqnarray*}
Thus, subtracting \eqref{3con1} from \eqref{3me1} implies
$$E_{i+1}\le\varphi_iE_i+h|\tau_i|,$$
where $E_i=|u_i-\textbf{u}_i|$ and
\begin{eqnarray*}
\varphi_i&=&1+h(w_1+w_2+w_3)L+(w_2a_{21}+w_3(a_{31}+a_{32}))h^2L^2+h^3L^3w_3a_{21}a_{32}+hL(w_2+(a_{21}w_2+a_{32}w_3)hL\\
&&+a_{21}a_{32}w_3h^2L^2)\frac{\varepsilon_{2i}^2(c_{2}h)^2}{2}   +hw_3L(1+hL(a_{31}+a_{32})+h^2L^2a_{21}a_{32})\frac{\varepsilon_{3i}^2(c_{3}h)^2}{2} \\ &&+h^2L^2(1+hLa_{21})a_{32}w_3\frac{\varepsilon_{2i}^2(c_{2}h)^2\varepsilon_{3i}^2(c_{3}h)^2}{4}.
\end{eqnarray*}
Let $$\mathcal{C}_3=\max_{i=0,\ldots,N-1}\left(\sup_{0<h\le(b-a)}\frac{2Lc_{2}^2h^2\mathcal{C}_{31}\varepsilon_{2i}^2+ 2Lc_{3}^2h^2\mathcal{C}_{32}\varepsilon_{3i}^2 +L^2h^5c_{2}^2c_{3}^2\mathcal{C}_{33}\varepsilon_{2i}^2\varepsilon_{3i}^2}{4e^{hL}}\right),$$
where
$\mathcal{C}_{31}=w_2+(c_{2}w_2+a_{32}w_3)hL+\frac{1}{6}h^2L^2$, $\mathcal{C}_{32}=1+hLw_3c_{3}+\frac{1}{6}h^2L^2$ and $\mathcal{C}_{33}=a_{32}w_3+\frac{1}{6}hL$. Since $\varepsilon_{2i}^2$ and $\varepsilon_{2i}^2$ are bounded for all $i=0,1,\ldots,N-1$, it follows that $\mathcal{C}_3<\infty$. Consequently, we have
\begin{eqnarray*}
\varphi_i &=&1+hL+\frac{1}{2}h^2L^2+\frac{1}{6}h^3L^3+hL\left(w_2+(c_{2}w_2+a_{32}w_3)hL+\frac{1}{6}h^2L^2\right)\frac{\varepsilon_{2i}^2(c_{2}h)^2}{2}   \\ &&
+hL\left(1+hLw_3c_{3}+\frac{1}{6}h^2L^2\right)\frac{\varepsilon_{3i}^2(c_{3}h)^2}{2} +h^2L^2\left(a_{32}w_3+\frac{1}{6}hL\right)\frac{\varepsilon_{2i}^2(c_{2}h)^2\varepsilon_{3i}^2(c_{3}h)^2}{4}\\
&\le&(1+\mathcal{C}_3h)e^{hL}\\
&\le&e^{(L+\mathcal{C}_3)h}.
\end{eqnarray*}
Thus, by arguments similar to those in Theorem \ref{th1}, we can show that $\textbf{u}_N\to u_N=u(b)$ as $h$ tends to zero, i.e., the MQ-RBF-based RK$3$ method \eqref{3me1} is convergent.
\begin{theorem}
Assume that $\varepsilon_{2i}^2$, $\varepsilon_{3i}^2$ and $\varepsilon_{4i}^2$ are bounded for all $i=0,1,\ldots,N-1$. Then the MQ-RBF-based RK$4$ method \eqref{4me1}, satisfying \eqref{pv41} and \eqref{pv42}, is convergent.
\end{theorem}
\textit{Proof.} For $t=t_i$, let us define
\begin{align*}
K_1(v)&=f(t_{i},v),\\
K_2(v)&=f\left(t_{i}+c_2 h, \left(1+\frac{\varepsilon_{2i}^2(c_{2}h)^2}{2}\right)\left(v+ha_{21}K_1(v)\right)\right),\\
K_3(v)&=f\left(t_{i}+c_3 h, \left(1+\frac{\varepsilon_{3i}^2(c_{3}h)^2}{2}\right)\big(v+h(a_{31}K_1(v)+a_{32}K_2(v))\big)\right)\\
K_4(v)&=f\left(t_{i}+c_4 h, \left(1+\frac{\varepsilon_{4i}^2(c_{4}h)^2}{2}\right)\big(v+h(a_{41}K_1(v)+a_{42}K_2(v)+a_{43}K_3(v))\big)\right),
\end{align*}
so that, using \eqref{e4}, we can write the exact solution of \eqref{ivp} as
\begin{equation}\label{4con1}
u_{i+1}=u_i+h(w_1K_1(u_i)+w_2K_2(u_i)+w_3K_3(u_i)+w_4K_4(u_i))+h\tau_i.
\end{equation}
Moreover, it follows from the Lipschitz continuity of $f(t,u)$ with respect to $u$ that
\begin{eqnarray*}
|K_1(u_i)-K_1(\textbf{u}_i)|&\le& L|u_i-\textbf{u}_i|,\\
|K_2(u_i)-K_2(\textbf{u}_i)| &\le& L(1+ha_{21}L)\left(1+\frac{\varepsilon_{2i}^2(c_{2}h)^2}{2}\right)|u_i-\textbf{u}_i|,\\
|K_3u(_i)-K_3(\textbf{u}_i)|&\le&L\left(1+\frac{\varepsilon_{3i}^2(c_{3}h)^2}{2}\right)\left(1+h(a_{31}+a_{32})L+a_{21}a_{32}h^2L^2+(1+hLa_{21})a_{32}hL\frac{\varepsilon_{2i}^2(c_{2}h)^2}{2}\right)|u_i-\textbf{u}_i|,
\end{eqnarray*}
and
\begin{eqnarray*}
|K_4(u_i)-K_4(\textbf{u}_i)|&\leq & L \left(\frac{1+\varepsilon_{4i}^{2}(c_{4}h)^2}{2}\right) \left[ 1+hL\left( a_{41}+a_{42}\left(\frac{1+\varepsilon_{2i}^{2}(c_{2}h)^2}{2}\right)+a_{43}\left(\frac{1+\varepsilon_{3i}^{2}(c_{3}h)^2}{2}\right) \right) \right.\\
&&+ h^{2}L^{2} \left( a_{42}a_{21}\left(\frac{1+\varepsilon_{2i}^{2}(c_{2}h)^2}{2}\right)+\left(\frac{1+\varepsilon_{3i}^{2}(c_{3}h)^2}{2}\right) \left( a_{43}a_{31}+a_{43}a_{32}\left(\frac{1+\varepsilon_{2i}^{2}(c_{2}h)^2}{2}\right)\right) \right) \\
&&\left.+h^{3}L^{3} \left( \left(\frac{1+\varepsilon_{3i}^{2}(c_{3}h)^2}{2}\right) 
\left( a_{43}a_{21}a_{31}+a_{43}a_{21}a_{32} \frac{\varepsilon_{2i}^{2}(c_{2}h)^{2}}{2} \right) \right) \right]|u_{i}-\textbf{u}_{i}|.
\end{eqnarray*}
Using \eqref{4con1} and \eqref{4me1} and employing the above inequalities, we can obtain
$$E_{i+1}=\varphi_iE_i+h|\tau_{i}|,$$
where $E_i=|u_i-\textbf{u}_i|$ and
$ \varphi_i\le 1+hL+\frac{1}{2}h^2L^2+\frac{1}{6}h^3L^3+\frac{1}{24}h^4L^4+\mathcal{C}_4he^{hL}\le e^{(L+\mathcal{C}_4)h},$
where, $$\mathcal{C}_4=\max_{i=0,\ldots,N-1}\left(\sup_{0<h\le(b-a)}\frac{4\mathcal{C}_{41}\varepsilon_{2i}^2+4\mathcal{C}_{42}\varepsilon_{3i}^2+4\mathcal{C}_{43}\varepsilon_{4i}^2+2\mathcal{C}_{44}\varepsilon_{2i}^2\varepsilon_{3i}^2+2\mathcal{C}_{45}\varepsilon_{2i}^2\varepsilon_{4i}^2+2\mathcal{C}_{46}\varepsilon_{3i}^2\varepsilon_{4i}^2+\mathcal{C}_{47}\varepsilon_{2i}^2\varepsilon_{3i}^2\varepsilon_{4i}^2}{8}\right),$$
with
\begin{eqnarray*}
\mathcal{C}_{41}&=&h^2Lc_2^2\left(w_2+(c_{2}w_2+a_{32}w_3+a_{42}w_4)hL+c_{2}(a_{32}w_3+a_{42}w_4)h^2L^2+\frac{1}{24}h^3L^3\right),\\    \mathcal{C}_{42}&=&h^2Lc_3^2(w_3+hL(w_3c_3+w_4a_{43})+h^2L^2(w_3c_{2}a_{32}+w_4a_{31}a_{43})+h^3L^3w_4a_{43}c_{2}a_{31}), \\ 
\mathcal{C}_{43}&=&h^2Lc_4^2(w_4+hLw_4c_4+h^2L^2w_4(c_{2}a_{42}+a_{43}c_3)+h^3L^3w_4a_{43}c_{2}a_{31}),\\
\mathcal{C}_{44}&=&h^3L^2c_2^2c_3^2(1+hLc_{2})a_{32}(w_3+w_4a_{43}),\\
\mathcal{C}_{45}&=&h^5L^2c_2^2c_4^2w_4a_{42}(1+c_{2}hL),\\
\mathcal{C}_{46}&=&h^5L^2c_3^2c_4^2w_4a_{43}(1+hLa_{31}),\\
\mathcal{C}_{47}&=&h^7L^2c_2^2c_3^2c_4^2w_4a_{43}a_{32}(1+hLc_{2}),
\end{eqnarray*}
is bounded under conditions \eqref{pv41}, \eqref{pv42} and boundedness of $\varepsilon_{2i}^2$, $\varepsilon_{3i}^2$ and $\varepsilon_{4i}^2$.
Therefore, $\textbf{u}_N \to u_N = u(b)$ as $h$ tends to zero. This implies that the MQ-RBF-based RK4 method converges. The proof is complete.
\section{Stability regions for MQ-RBF-based RK Methods}\label{sec4}
\noindent
In Section \ref{sec3}, we have demonstrated the convergence of MQ-RBF-based RK methods with up to four stages over a bounded interval when the step size $h$ approaches zero. However, in practice, convergence alone does not guarantee accurate results for a given step size $h > 0$. It is significantly influenced by the stability region of the numerical method. 
In this section, we graphically illustrate the stability regions of all MQ-RBF-based RK methods with up to four stages and compare them with those of the corresponding classical RK method. 

Any one-step method, when applied to the test problem $u'=\lambda u$, can be expressed in the form 
$$\textbf{u}_{i+1}=\mathcal{R}(z)\textbf{u}_i,$$ 
where $z=\lambda h$ and $\mathcal{R}(z)$ is the stability function.
For a consistent one-step method, the stability function $\mathcal{R}(z)$ must approximate $e^z$ as $z \to 0$. Furthermore, if the method is of order $p$, then
\begin{eqnarray*}
\mathcal{R}(z)-e^z=\mathcal{O}(z^{p+1}) \ \ \text{as} \ \ h\to 0.
\end{eqnarray*}
It is well known that classical RK methods are one-step methods, and the stability functions for the RK2, RK3, and RK4 methods are, respectively, given as 
\begin{eqnarray*}
\mathcal{R}(z)&=&1+z+\frac{z^2}{2},\\
\mathcal{R}(z)&=&1+z+\frac{z^2}{2}+\frac{z^3}{6},\\
\mathcal{R}(z)&=&1+z+\frac{z^2}{2}+\frac{z^3}{6}+\frac{z^4}{24}.\\
\end{eqnarray*}
It is clear that all MQ-RBF-based RK methods developed up to four stages are one-step methods. We now present the stability functions of the MQ-RBF-based RK methods.

\subsection{MQ-RBF-based RK2 method}
\noindent
When the MQ-RBF-based RK2 method presented in Section \ref{MR1} is applied to the test problem $u'=\lambda u$, the optimized value of the shape parameter $\varepsilon_{2i}^2$ appearing in the method is given by $\frac{z^2}{h^2}$. Thus, the stability function for the MQ-RBF-based RK2 method can be obtained as
$$\mathcal{R}(z)
=1+z+\frac{z^2}{2}+\frac{z^3}{6}+\frac{z^4}{9}$$
It is clear that
\begin{eqnarray*}
\mathcal{R}(z)-e^z=\mathcal{O}(z^4),
\end{eqnarray*}
which confirms the third-order convergence of the method.

\subsection{MQ-RBF-based RK3 method}
\noindent
When the MQ-RBF-based RK3 method presented in Section~\ref{MR2} is applied to the test problem $u'=\lambda u$, the optimized values of the shape parameters are obtained as listed in Table \ref{tab1}.

\begin{table}[h!]
\caption{Optimized values of the shape parameters when the MQ-RBF-based RK3 method is applied to the test problem $u'=\lambda u$.}
\centering
\begin{tabular}{|c|ccc|}
\hline\hline
\textbf{MQ-RBF-based RK3 methods}& 
\multicolumn{1}{c}{\textbf{Shape parameters}}&&\\
\cline{2-3}
(Section \ref{MR2})& $\varepsilon_{2i}^2$ &$\varepsilon_{3i}^2$&\\
\hline
\textbf{B1}  \hspace{0.5cm} &---&---&\\
\textbf{B2 (a)}  & $\frac{(7-\sqrt{33})z^2}{4h^2}$&$-\frac{z^2}{h^2}$&\\
\textbf{B2 (b)}  & $\frac{(7+\sqrt{33})z^2}{4h^2}$&$-\frac{z^2}{h^2}$&\\
\textbf{B3 (a)} &$\frac{z^2}{h^2}$&$-\frac{z^2}{h^2}$&\\
\textbf{B3 (b)}  &$\frac{3z^2}{h^2}$&$-\frac{3z^2}{5h^2}$&\\
\textbf{B4}  \hspace{0.5cm} &$\frac{4z^2}{3h^2}$&$-\frac{4z^2}{9h^2}$&\\
\hline\hline
\end{tabular}
\label{tab1}
\end{table}
Using each set of parameter values from Section~\ref{MR2}, together with the corresponding shape parameters from Table~\ref{tab1}, the stability function for the MQ-RBF-based RK3 method can be obtained from \eqref{2me2} as follows:
\begin{description}
\item[B1] \begin{eqnarray*}
\mathcal{R}(z)&=&1+z+\frac{z^{2}}{2}+\frac{z^{3}}{6}-\frac{z^{5}}{32}-\frac{z^{6}}{192}-\frac{z^{7}}{134}.
\end{eqnarray*}
\item[B2] 
\begin{description}
\item[(a)] \begin{eqnarray*}
	\mathcal{R}(z)&=&1+z+\frac{z^{2}}{2}+\frac{z^{3}}{6}+\frac{z^{4}}{24}-\left(\frac{1}{128}-\frac{\sqrt{33}}{384}\right)z^{5}-\left(\frac{23}{1152}-\frac{11\sqrt{33}}{3456}\right)z^{6}-\left(\frac{7}{864}-\frac{\sqrt{33}}{864}\right)z^7.
\end{eqnarray*}
\item[(b)] 
\begin{eqnarray*}
	\mathcal{R}(z)&=&1+z+\frac{z^{2}}{2}+\frac{z^{3}}{6}+\frac{z^{4}}{24}-\left(\frac{1}{128}+\frac{\sqrt{33}}{384}\right)z^{5}-\left(\frac{23}{1152}+\frac{11\sqrt{33}}{3456}\right)z^{6}-\left(\frac{7}{864}+\frac{\sqrt{33}}{864}\right)z^7.
\end{eqnarray*}
\end{description}
\item[B3] 
\begin{description}
\item[(a)] \begin{eqnarray*}
	\mathcal{R}(z)
	&=&1+z+\frac{z^{2}}{2}+\frac{z^{3}}{6}+\frac{z^{4}}{24}+\frac{z^{5}}{48}-\frac{z^{6}}{864}-\frac{z^{7}}{864}.
\end{eqnarray*}
\item[(b)] \begin{eqnarray*}
	\mathcal{R}(z)
	&=&1+z+\frac{z^{2}}{2}+\frac{z^{3}}{6}+\frac{z^{4}}{24}-\frac{z^{5}}{144}-\frac{5z^{6}}{288}-\frac{5z^{7}}{864}.
\end{eqnarray*}
\end{description}
\item[B4] \begin{eqnarray*}
\mathcal{R}(z)
&=&1+z+\frac{z^{2}}{2}+\frac{z^{3}}{6}+\frac{z^{4}}{24}+\frac{z^{5}}{144}-\frac{z^{6}}{144}-\frac{z^{7}}{288}.
\end{eqnarray*}
\end{description}
It is evident from the expressions of $\mathcal{R}(z)$ that for each case except \textbf{B1}, we have
\begin{eqnarray*}
\mathcal{R}(z)-e^z=\mathcal{O}(z^5),
\end{eqnarray*}
which confirms fourth-order convergence. However, for any choice of shape parameters, the case \textbf{B1} yields 
\begin{eqnarray*}
\mathcal{R}(z)-e^z=\mathcal{O}(z^4),
\end{eqnarray*}
indicating that, for the given test problem, the MQ-RBF-based RK3 method with parameter values from \textbf{B1} (see Section~\ref{MR2}) achieves only third-order accuracy.

\subsection{MQ-RBF-based RK4 method}
\noindent
For $f=\lambda u$, the optimized values of the shape parameters for the MQ-RBF-based RK3 method presented in Section~\ref{MR3} are listed in Table~\ref{tab2}. 
\begin{table}[h!]
\caption{Optimized values of the shape parameters when the MQ-RBF-based RK4 method is applied to the test problem $u'=\lambda u$.}
\centering
\begin{tabular}{|c|cccc|}
\hline\hline
\textbf{MQ-RBF-based RK4 methods}& 
\multicolumn{1}{c}{\textbf{Shape parameters}}&&&\\
\cline{2-4}
(Section \ref{MR3})& $\varepsilon_{2i}^2$ &$\varepsilon_{3i}^2$&\hspace{1cm}$\varepsilon_{4i}^2$&\\
\hline
\textbf{C1}  &$-\frac{55z^2}{12h^2}$&$\frac{55z^2}{18h^2}$&\hspace{1cm}$-\frac{5z^2}{6h^2}$&\\
\textbf{C2 }  & $\frac{8z^2}{3h^2}$&$-\frac{4z^2}{9h^2}$&\hspace{1cm}$\frac{4z^2}{15h^2}$&\\
\hline\hline
\end{tabular}
\label{tab2}
\end{table}

Using the values of the shape parameters from Table~\ref{tab2} along with the corresponding set of parameter values from Section~\ref{MR2}, the stability function for the MQ-RBF-based RK4 method can be obtained from \eqref{3me2} as follows:
\begin{description}
\item[C1]
$$R(z)=1+z+\frac{z^{2}}{2}+\frac{z^{3}}{6}+\frac{z^{4}}{24}+\frac{z^{5}}{120}-\frac{1763z^{6}}{17280}-\frac{209z^{7}}{4320}-\frac{1001z^{8}}{86400}-\frac{121z^{9}}{13824}-\frac{121z^{10}}{34560}.$$
\item[C2] 
$$R(z)=1+z+\frac{z^{2}}{2}+\frac{z^{3}}{6}+\frac{z^{4}}{24}+\frac{z^{5}}{120}-\frac{37z^{6}}{21600}-\frac{z^{7}}{540}-\frac{7z^{8}}{27000}-\frac{z^{9}}{6750}-\frac{z^{10}}{27000}. $$
\end{description}
Clearly, for both cases, we have
\begin{eqnarray*}
\mathcal{R}(z)-e^z=\mathcal{O}(z^6),
\end{eqnarray*}
thus verifying the fifth-order accuracy of the scheme.

\begin{figure}[h!]
\centering
\includegraphics[width=0.45\linewidth]{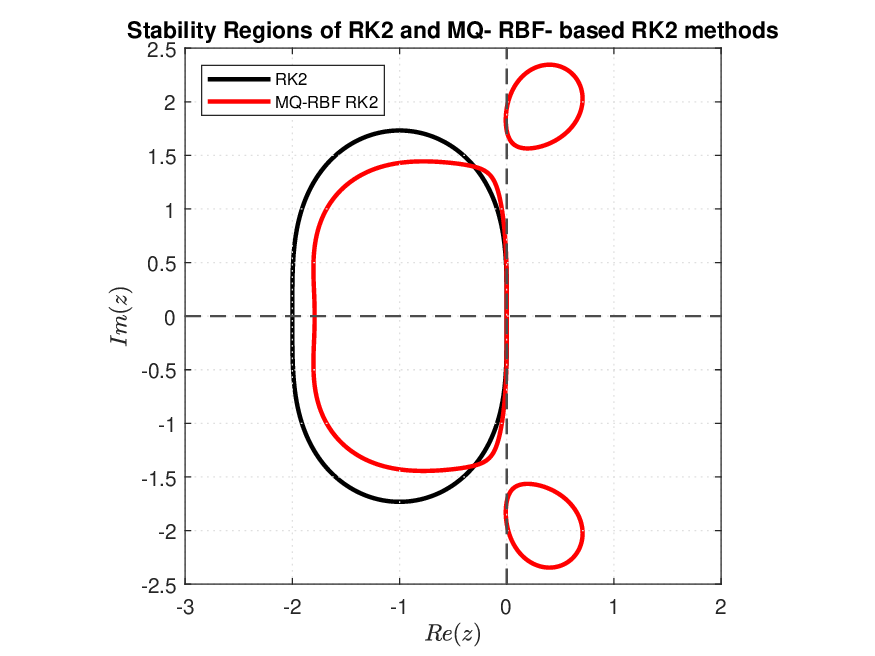} 
\includegraphics[width=0.45\linewidth]{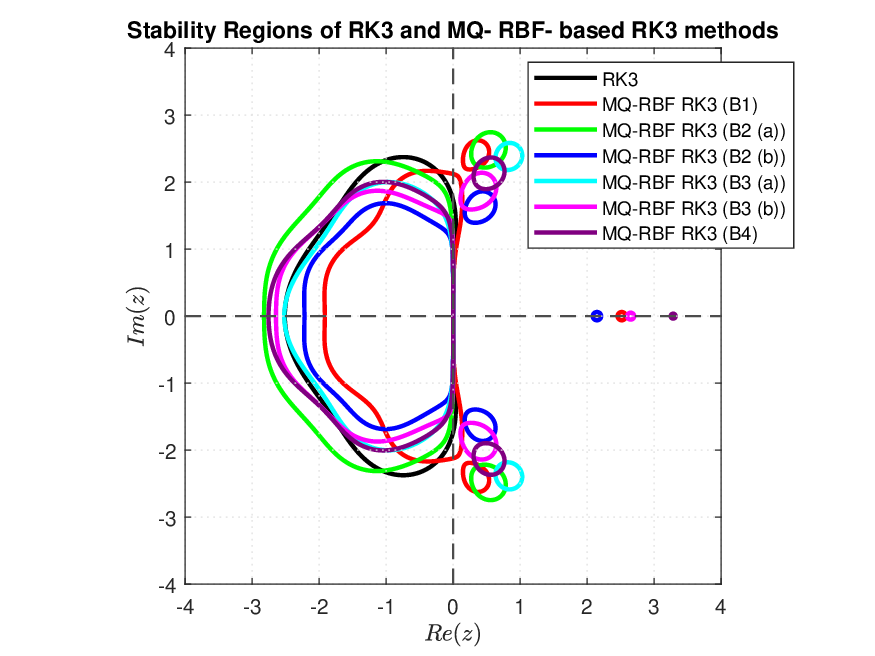}
\includegraphics[width=0.45\linewidth]{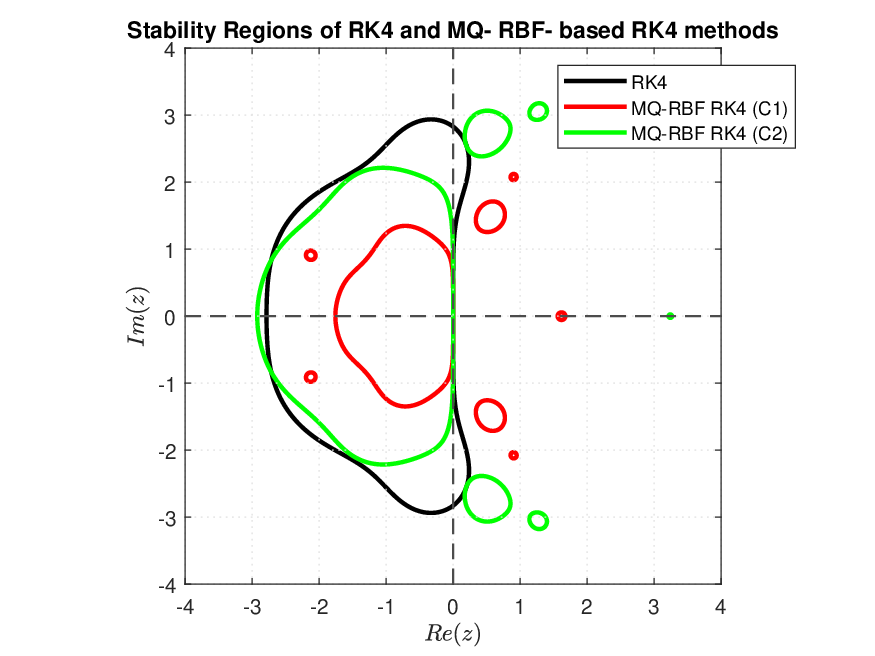}
\caption{Stability regions of classical RK and MQ-RBF-based RK methods: two-stage methods (top left), three-stage method (top right) and four-stage methods (bottom center).}
\label{fig1}
\end{figure}
Next, we plot the stability regions of the classical RK method and the corresponding MQ-RBF-based RK method to perform a comparative study of their stability intervals. Here, the term stability interval refers to the portion of the negative real axis that lies within the stability region. From the top-left plot of Figure~\ref{fig1}, we observe that the stability interval of the MQ-RBF-based RK2 method is slightly smaller than that of the classical RK2 method. The stability regions for the classical RK3 method and various MQ-RBF-based RK3 methods, corresponding to different parameter values given in Section~\ref{MR2}, are shown in the top right plot of Figure~\ref{fig1}. From the largest to the smallest stability interval, the corresponding regions are associated with MQ-RBF--RK3 (\textbf{B2(a)}), MQ-RBF--RK3 (\textbf{B4}), MQ-RBF--RK3 (\textbf{B3(b)}), RK3, MQ-RBF--RK3 (B3(a)), MQ-RBF--RK3 (\textbf{B2(b)}), and MQ-RBF--RK3 (\textbf{B1}). The bottom center plot of Figure~\ref{fig1} shows the stability regions for the classical RK4 method and the MQ-RBF-based RK4 methods corresponding to different parameter values given in Section~\ref{MR3}. From the largest to the smallest stability interval, the corresponding regions are associated with MQ-RBF--RK4 (\textbf{C2}), RK4, and MQ-RBF--RK4 (\textbf{C1}), respectively. 
\section{Numerical Results}\label{sec5}
In this section, we present numerical examples based on initial value problems to demonstrate the improved precision of the proposed multiquadratic RBF-based Runge–Kutta (MQ-RBF--RK) methods.We have evaluated the numerical error and the order of convergence of the proposed methods and compared the results with those of the standard Runge–Kutta method. All computations are performed with MATLAB, version R2024b.\\

\begin{example} \label{6.1}
In this example, we assume $f$ to be a function independent of $t$. Consider the following IVP:
\begin{equation*}
\frac{du}{dt}=-u^2, ~~~t\geq0,
\end{equation*}
with the initial condition $u(0)=1$. The exact solution for this IVP is given by $u(t)=\frac{1}{t+1}$. To compute the numerical solution, we consider the final time as $T=1$.

\begin{table}[h!]
\caption{Optimal value of the shape parameter $\varepsilon_{2i}^2$ when the MQ-RBF--RK methods are applied to Example~\ref{6.1}.}
\centering
\begin{tabular}{|c|ccccccccc|}
	\hline\hline
	\textbf{ Method }  &\textbf{MQ-RBF--RK2}&&\multicolumn{6}{c}{\textbf{MQ-RBF--RK3}}&\\
	\cline{4-9}
	&&&(\textbf{B1})&(\textbf{B2 (a)})&(\textbf{B2 (b)})&(\textbf{B3 (a)})&(\textbf{B3 (b)})&\hspace{0.5cm}(\textbf{B4})&\\
	$\varepsilon_{2i}^2$  & $2u^2$&&$450u^2 \ ${\footnote{For this particular problem, using the parameters of the MQ-RBF--RK3 (\textbf{B1}) method, the coefficient of $\varepsilon_{2i}^2$ in \eqref{e3} becomes zero. As a result, a specific value for the method's shape parameter cannot be determined, and we have therefore chosen it arbitrarily.}}&$\left(\frac{13}{4}-\frac{5\sqrt{33}}{12}\right)u^2$&$\left(\frac{13}{4}+\frac{5\sqrt{33}}{12}\right)u^2$&$\frac{2}{3}u^2$&$6u^2$&\hspace{0.5cm}$\frac{8}{3}u^2$&\\
	\hline \hline
	\textbf{Method}   &&\multicolumn{4}{c}{\textbf{MQ-RBF--RK4}}&&&&  \\
	\cline{3-6}
	&& (\textbf{C1}$+$)&(\textbf{C1}$-$)&(\textbf{C2$+$})&(\textbf{C2$-$})&&&&\\
	$\varepsilon_{2i}^2$ && $\frac{34+\sqrt{2066}}{14}u^2$&$\frac{34-\sqrt{2066}}{14}u^2$&$(-4+2\sqrt{23})u^2$&$(-4-2\sqrt{23})u^2$&&&&\\
	\hline \hline
\end{tabular}
\label{tab3}
\end{table}

The values of the shape parameter $(\varepsilon_{2i}^2)$ for each of the MQ-RBF-based RK methods are listed in Table~\ref{tab3}. All other shape parameters follow the relationship provided in Section~\ref{sec2}. Using these values, we have calculated the global error and the order of convergence of each method at the final time $T$ for various values of the number of steps $N$ (see Tables \ref{Eg1rk2},\ref{Eg1rk3} and \ref{Eg1rk4}). For comparison, we have also provided the error and order of convergence of the classical RK methods, which emerges as the limiting case of the corresponding MQ-RBF–RK methods when the RBF shape parameter approaches zero. From Tables \ref{Eg1rk2}, \ref{Eg1rk3}, and \ref{Eg1rk4}, it is evident that for each $N$, the classical RK method with $s$‑stages achieves its usual order of convergence $s$, while our proposed $s$‑stage MQ-RBF-based RK methods attain an improved order of convergence $s + 1$. Also, except for the MQ-RBF--RK3 (\textbf{B1}) method, and for $N=20$, the MQ-RBF--RK4 (\textbf{C2$-$}) method, the global errors in our MQ-RBF-based RK methods are significantly lower for each $N$ than those observed in the corresponding RK methods. In addition, the MQ-RBF--RK2 method in Table~\ref{Eg1rk2} exhibits lower errors than most of the RK3 methods in Table~\ref{Eg1rk3}. Thus, we can conclude that the MQ-RBF--RK2 method is more accurate than the RK3 methods.
\begin{table}[h!]
\caption{Global error and order of convergence of RK2 and MQ-RBF-based RK2 methods applied to Example~\ref{6.1}.} 
\begin{center}
	\begin{tabular}{ccccccccccccccc}
		\hline\hline
		&&& \multicolumn{2}{c}{\textbf{RK2}}  &&& \multicolumn{2}{c}{\textbf{MQ-RBF--RK2}} &&& \\
		\cline{3-5}\cline{7-9} 
		$N$    && &     error    &    order &&&   error  &    order   \\
		\hline\hline			
		20  && & 2.20e-04  & .... &&& 1.21e-06 & .... \\
		40  & && 5.36e-05  & 2.0410 &&& 1.58e-07 & 2.9429 \\
		80 & && 1.32e-05   & 2.0204 &&& 2.00e-08 & 2.9754\\
		160 & && 3.28e-06  & 2.0102 &&& 2.52e-09& 2.9886\\
		320 & && 8.17e-07  & 2.0051 &&&3.17e-10 & 2.9945\\
		\hline\hline
	\end{tabular}
\end{center}
\label{Eg1rk2}
\end{table}
\begin{table}[h!]
\caption{Global error and order of convergence of RK3 and MQ-RBF-based RK3 methods applied to Example~\ref{6.1}.}
\begin{center}
	\begin{tabular}{cccccccccccccccccccccc}
		\hline\hline
		&&& \multicolumn{2}{c}{\textbf{RK3 (B1)}}  &&& \multicolumn{2}{c}{\textbf{MQ-RBF--RK3 (B1)}} &&& \multicolumn{2}{c}{\textbf{RK3 (B2(a))}}& &&\multicolumn{2}{c}{\textbf{MQ-RBF--RK3 (B2(a))}}\\
		\cline{4-5}\cline{8-9}\cline{12-13}\cline{16-17}
		$N$ &&& error & order &&& error & order &&&  error & order &&&error & order\\
		\hline\hline
		20  &&& 2.16e-06 & ... &&& 5.47e-03 & ...   &&& 5.76e-06  &  ... &&& 1.19e-07&...\\
		40  &&& 2.57e-07 & 3.0752 &&&  3.26e-04 &  4.0696 &&&6.93e-07 &  3.0558&&&7.19e-09 &4.0524\\
		80  &&& 3.13e-08 & 3.0363 &&& 1.94e-05  & 4.0677  &&& 8.49e-08 &  3.0280&&&4.41e-10 &4.0265\\
		160 &&& 3.86e-09 & 3.0178 &&& 1.19e-06  &  4.0348 &&& 1.05e-08 &  3.0140&&& 2.73e-11&4.0133\\
		320 &&& 4.80e-10 & 3.0089 &&&  7.33e-08 &  4.0141 &&& 1.31e-09 & 3.0070&&& 1.70e-12&4.0069\\
		\hline\hline
		&&& \multicolumn{2}{c}{\textbf{RK3 (B2(b))}}  &&& \multicolumn{2}{c}{\textbf{MQ-RBF--RK3 (B2(b))}} &&& \multicolumn{2}{c}{\textbf{RK3 (B3(a))}}& &&\multicolumn{2}{c}{\textbf{MQ-RBF--RK3 (B3(a))}}\\
		\cline{4-5}\cline{8-9}\cline{12-13}\cline{16-17}
		$N$ &&& error & order &&& error & Order &&&  error & order &&&error & order\\
		\hline\hline
		20  &&& 3.68e-06 & ... &&& 1.31e-07& ...   &&& 4.14e-06 &  ... &&& 5.48e-08&...\\
		40  &&& 4.46e-07 & 3.0450 &&& 8.14e-09 & 4.0089 &&& 5.03e-07 &  3.0410&&& 3.36e-09&4.0266\\
		80  &&& 5.49e-08 & 3.0225 &&& 5.07e-10 & 4.0062 &&& 6.19e-08 &  3.0208&&& 2.08e-10&4.0139\\
		160 &&& 6.81e-09& 3.0112 &&& 3.16e-11 & 4.0035 &&& 7.69e-09 &  3.0105&&&1.29e-11 &4.0071\\
		320 &&& 8.48e-10 & 3.0056 &&& 1.97e-12 & 4.0019 &&& 9.57e-10 & 3.0052&&& 8.07e-13&4.0043\\
		\hline\hline
		&&& \multicolumn{2}{c}{\textbf{RK3 (B3(b))}}  &&& \multicolumn{2}{c}{\textbf{MQ-RBF--RK3 (B3(b))}} &&& \multicolumn{2}{c}{\textbf{RK3 (B4)}}& &&\multicolumn{2}{c}{\textbf{MQ-RBF--RK3 (B4)}}\\
		\cline{4-5}\cline{8-9}\cline{12-13}\cline{16-17}
		$N$ &&& error & Order &&& error & order &&&  error & order &&&error & order\\
		\hline\hline
		20  &&& 4.16e-06 & ... &&& 1.21e-07& ...   &&& 4.16e-06 &  ... &&& 8.87e-08&...\\
		40  &&& 5.04e-07 & 3.0450 &&& 7.4e-09 & 4.0257 &&& 5.04e-07 &  3.0453&&&5.41e-09 &4.0354\\
		80  &&& 6.20e-08 & 3.0226 &&& 4.58e-10 & 4.0140 &&& 6.20e-08 &  3.0227&&& 3.34e-10&4.0183\\
		160 &&& 7.69e-09 & 3.0113 &&& 2.85e-11 & 4.0073 &&& 7.69e-09 &  3.0114&&& 2.07e-11&4.0093\\
		320 &&& 9.57e-10 & 3.0056 &&& 1.78e-12 & 4.0036 &&& 9.58e-10 & 3.0057&&&1.29e-12 &4.0048\\
		\hline \hline
	\end{tabular}
\end{center}
\label{Eg1rk3}
\end{table}
\begin{table}[h!]
\caption{Global error and order of convergence of RK4 and MQ-RBF-based RK4 methods applied to Example~\ref{6.1}.}
\begin{center}
	\begin{tabular}{ccccccccccccccccccc}
		\hline\hline
		&&& \multicolumn{2}{c}{\textbf{RK4 (C1)}}  &&& \multicolumn{2}{c}{\textbf{MQ-RBF--RK4 (C1$+$)}} &&& \multicolumn{2}{c}{\textbf{MQ-RBF--RK4 (C1$-$)}}& &&\\
		\cline{4-5}\cline{8-9}\cline{12-13}
		$N$ &&& error & order &&& error & order &&&  error & order &&&\\
		\hline\hline
		20  &&& 1.69e-08 & ... &&& 1.21e-08& ...   &&& 5.08e-09 &  ... \\
		40  &&& 1.09e-09 & 3.9501 &&& 3.55e-10 & 5.0906 &&& 1.49e-10 &  5.0912\\
		80  &&& 6.93e-11 & 3.9796 &&& 1.07e-11 & 5.0455 &&& 4.51e-12 &  5.0457\\
		160 &&& 4.35e-12 & 3.9909 &&& 3.30e-13 & 5.0225 &&& 1.38e-13 &  5.0274\\
		320 &&& 2.73e-13 & 3.9951 &&& 1.03e-14 & 5.0073 &&& 4.55e-15 & 4.9250\\
		\hline\hline
		&&& \multicolumn{2}{c}{\textbf{RK4 (C2)}}  &&& \multicolumn{2}{c}{\textbf{MQ-RBF--RK4 (C2$+$)}} &&& \multicolumn{2}{c}{\textbf{MQ-RBF--RK4 (C2$-$)}} \\
		\cline{4-5}\cline{8-9} \cline{12-13}
		$N$ &&& error & order &&& error & order &&&  error & order \\
		\hline\hline
		20  &&& 3.74e-08 & ... &&&2.03e-09 & ...   &&& 4.77e-08 &  ... \\
		40  &&& 2.30e-09 & 4.0237 &&& 5.97e-11 & 5.0874 &&& 1.37e-09 &  5.1183\\
		80  &&& 1.42e-10 & 4.0132 &&& 1.81e-12 & 5.0438 &&& 4.12e-11 &  5.0593\\
		160 &&& 8.85e-12 & 4.0069 &&& 5.57e-14& 5.0232 &&& 4.26e-12 &  5.0291\\
		320 &&& 5.52e-13 & 4.0038 &&& 1.67e-15 & 5.0632 &&& 3.91e-14 & 5.0110\\ 
		\hline\hline
	\end{tabular}
\end{center}
\label{Eg1rk4}
\end{table}

\end{example}
\begin{example}\label{6.2}
Next, we consider a stiff problem given as
\begin{equation*}
\frac{du}{dt}=-4t^3u^2, ~~~t\geq-10,
\end{equation*}
with the initial condition $u(-10)=\frac{1}{10001}$. The exact solution for this IVP is given by $u(t)=\frac{1}{t^4+1}$. We approximate the numerical solution at the final time $T=0$ and compute the global error and the order of convergence for each proposed method using different values of the number of steps $N$.

\begin{table}[h!]
\caption{Optimal value of the shape parameter $\varepsilon_{2i}^2$ when the MQ-RBF--RK methods are applied to Example~\ref{6.2}.}
\centering
\begin{tabular}{l cccc}
	\hline\hline
	\textbf{Methods}& &&
	$\varepsilon_{2i}^2$\\
	\hline
	MQ-RBF--RK2 \hspace{0.5cm}& &&$32t^6u^2-12t^2u$\\
	MQ-RBF--RK3 (\textbf{B1})  \hspace{0.5cm} &&&$\frac{128}{3}t^{10}u^3+16t^6u^2-12t^2u$\\
	MQ-RBF--RK3 (\textbf{B2 (a)}) &&& $\frac{24t^2u((2048172224074424 t^8 u^2-973179214634506 t^4 u + 205114630606597))}{3587275474121788 t^4 u - 615343891819791}$\\
	MQ-RBF--RK3 (\textbf{B2 (b)})  &&& $\frac{24t^2 u(1329527496453448t^8u^2 - 434195668918774 t^4 u - 64377142251269)}{353374199827396 t^4 u + 193131426753807}$\\
	MQ-RBF--RK3 (\textbf{B3 (a)}) &&&$\frac{(256 t^{10} u^3 - 48 t^6 u^2 + 12 t^2 u)}{(24 t^4 u - 3)}$\\
	MQ-RBF--RK3 (\textbf{B3 (b)})  &&&$\frac{768 t^{10} u^3 - 336 t^6 u^2 - 12 t^2 u}{8 t^4 u + 3}$\\
	MQ-RBF--RK3 (\textbf{B4})  \hspace{0.5cm} &&&$\frac{1024 t^{12} u^3 - 432 t^8*u^2 + 24t^4 u - 1}{24t^6 u + 3 t^2}$\\
	MQ-RBF--RK4 (\textbf{C1$+$})  \hspace{0.5cm} &&&$ \frac{1088u^2t^8 + \sqrt{2115584t^16u^4 - 1772544t^{12}u^3 + 922656t^8u^2 + 79200t^4u + 1089} - 1200t^4u - 33}{28t^2} $\\
	MQ-RBF--RK4 (\textbf{C1$-$})  \hspace{0.5cm} &&&$ \frac{1088u^2t^8 - \sqrt{2115584t^16u^4 - 1772544t^{12}u^3 + 922656t^8u^2 + 79200t^4u + 1089} - 1200t^4u - 33}{28t^2} $\\
	MQ-RBF--RK4 (\textbf{C2$+$})  \hspace{0.5cm} &&&$ \frac{-64u^2t^8 + 2\sqrt(5888t^16u^4 + 4416t^{12}u^3 + 2064t^8u^2 + 468t^4u + 9) - 120t^4u - 6}{t^2} $\\
	MQ-RBF--RK4 (\textbf{C2$-$})  \hspace{0.5cm} 
	&&&$ \frac{-64u^2t^8 - 2\sqrt(5888t^16u^4 + 4416t^{12}u^3 + 2064t^8u^2 + 468t^4u + 9) - 120t^4u - 6}{t^2} $\\
	\hline\hline
\end{tabular}
\label{tabs2}
\end{table}

The value of the shape parameter $\varepsilon_{2i}^2$ for each MQ-RBF--RK method applied to this example is provided in Table~\ref{tabs2}. From tables \ref{Eg2rk2}, \ref{Eg2rk3}, and \ref{Eg2rk4}, it is evident that each MQ-RBF-based RK method achieves an order of convergence one higher than that of the corresponding RK method and yields a lower global error. Specifically, an $s$-stage MQ-RBF-based RK method exhibits an order of convergence $s+1$ (approx.). 
However, certain methods show a degradation in accuracy for specific step sizes ($T/N$). For instance, with $N=400$, the MQ-RBF--RK3 (\textbf{B2(a)}) method exhibits an order $-1.4821$, followed by a return to the expected convergence rate. Similarly, the MQ-RBF--RK3 (\textbf{C1$+$}) and MQ-RBF--RK3 (\textbf{C2$+$}) methods initially demonstrate a consistent convergence rate of approximately $5$, which later degrades for $N=6400$. These temporary losses of convergence, observed in explicit methods with bounded stability regions, are likely due to the stiffness of the problem. Despite this, the global errors for these methods are lower than those produced by the corresponding classical RK methods. On a positive note, we also observe slight superconvergence in the MQ-RBF--RK3 (\textbf{C2$-$}) method for $N=6400$. Additionally, most of the $s$-stage MQ-RBF--RK methods yield a higher rate of convergence and better accuracy compared to the classical RK methods with $s+1$ stages.

\begin{table}[h!]
\caption{Global error and order of convergence of RK2 and MQ-RBF-based RK2 methods applied to Example~\ref{6.2}.}
\begin{center}
	\begin{tabular}{ccccccccccccccc}
		\hline\hline
		&&& \multicolumn{2}{c}{\textbf{RK2}}  &&& \multicolumn{2}{c}{\textbf{MQ-RBF--RK2}} &&& \\
		\cline{4-5}\cline{8-9} 
		$N$    && &     error    &    order  &&&   error  &    order   \\
		\hline\hline			
		200  && & 7.51e-01  & .... &&& 3.21e-02 & .... \\
		400  & && 4.40e-01  & 0.7726 &&& 4.10e-03 & 2.9725 \\
		800 & &&  1.67e-01  & 1.4016 &&& 5.22e-04 & 2.9710\\
		1600 & && 4.79e-02 & 1.7972 &&& 6.60e-05& 2.9835\\
		3200 & &&  1.25e-02  & 1.9425 &&&8.30e-06 & 2.9915\\
		6400& && 3.15e-03  & 1.9842  &&&1.04e-06 & 2.9957\\
		\hline\hline
	\end{tabular}
\end{center}
\label{Eg2rk2}
\end{table}
\begin{table}[h!]
\caption{Global error and order of convergence of RK3 and MQ-RBF-based RK3 methods applied to Example~\ref{6.2}.}
\begin{center}
	\begin{tabular}{cccccccccccccccccccccc}
		\hline\hline
		&&& \multicolumn{2}{c}{\textbf{RK3 (B1)}}  &&& \multicolumn{2}{c}{\textbf{MQ-RBF--RK3 (B1)}} &&& \multicolumn{2}{c}{\textbf{RK3 (B2(a))}}& &&\multicolumn{2}{c}{\textbf{MQ-RBF--RK3 (B2(a))}}\\
		\cline{4-5}\cline{8-9}\cline{12-13}\cline{16-17}
		$N$ &&& error & Order &&& error & Order &&&  error & Order &&&error & Order\\
		\hline\hline
		200  &&& 4.34e-02 & ... &&& 3.00e-03 & ...   &&& 6.61e-02 & ... &&& 1.32e-04 & ...\\
		400  &&& 5.85e-03 & 2.8907 &&&  1.88e-04 &  3.9962 &&&9.10e-03 & 2.85976 &&&  3.69e-04 &  -1.4821\\
		800  &&& 7.49e-04 & 2.9659 &&& 1.17e-05  & 4.0015  &&& 1.17e-03 & 2.9619 &&& 1.27e-05  & 4.8645\\
		1600 &&& 9.46e-05 & 2.9856 &&& 7.33e-07  &  4.0012 &&& 1.48e-04 & 2.9854 &&& 4.33e-07  &  4.8705\\
		3200 &&& 1.19e-05 & 2.9931 &&&  4.58e-08 &  4.0005 &&& 1.85e-05 & 2.9932 &&&  1.57e-08 &  4.7893\\
		6400 &&& 1.49e-06 & 2.9966 &&& 2.88e-09 & 3.9898 &&& 2.32e-06 & 2.9967 &&& 6.70e-10 & 4.5454\\ 
		\hline\hline
		&&& \multicolumn{2}{c}{\textbf{RK3 (B2(b))}}  &&& \multicolumn{2}{c}{\textbf{MQ-RBF--RK3 (B2(b))}} &&& \multicolumn{2}{c}{\textbf{RK3 (B3(a))}}& &&\multicolumn{2}{c}{\textbf{MQ-RBF--RK3 (B3(a))}}\\
		\cline{4-5}\cline{8-9}\cline{12-13}\cline{16-17}
		$N$ &&& error & Order &&& error & Order &&&  error & Order &&&error & Order\\
		\hline\hline
		200  &&& 4.41e-02  &  ... &&& 4.17e-04&...  &&& 6.75e-02 &  ... &&& 7.54e-04&...\\
		400  &&& 5.91e-03 &  2.8973&&&2.61e-05 &3.9988 &&& 9.27e-03 &  2.8644&&& 5.03e-05&3.9061\\
		800  &&& 7.55e-04 &  2.9701&&&1.63e-06 &4.0007 &&& 1.19e-03 &  2.9653&&& 3.37e-06&3.8988\\
		1600 &&& 9.51e-05 &  2.9879&&& 1.02e-07&4.0007 &&& 1.50e-04 &  2.9872&&&2.11e-07 &3.9954\\
		3200 &&& 1.19e-05 & 2.9943&&& 6.37e-09&3.9988 &&& 1.88e-05 & 2.9942&&& 1.39e-08&3.9268\\
		6400 &&& 1.50e-06 & 2.9972&&& 3.92e-10&4.0215 &&& 2.35e-06 & 2.9972&&& 8.78e-10&3.9848\\ 
		\hline\hline
		&&& \multicolumn{2}{c}{\textbf{RK3 (B3(b))}}  &&& \multicolumn{2}{c}{\textbf{MQ-RBF--RK3 (B3(b))}} &&& \multicolumn{2}{c}{\textbf{RK3 (B4)}}& &&\multicolumn{2}{c}{\textbf{MQ-RBF--RK3 (B4)}}\\
		\cline{4-5}\cline{8-9}\cline{12-13}\cline{16-17}
		$N$ &&& error & Order &&& error & Order &&&  error & Order &&&error & Order\\
		\hline\hline
		200  &&& 4.34e-02 & ... &&& 2.53e-04& ...   &&& 4.16e-06 &  ... &&& 9.60e-05&...\\
		400  &&& 5.83e-03 & 2.8985 &&& 1.58e-05 & 4.0065 &&& 5.04e-07 &  3.0453&&&7.00e-06 &3.7783\\
		800  &&& 7.43e-04 & 2.9704 &&& 9.82e-07 & 4.0049 &&& 6.20e-08 &  3.0227&&& 4.71e-07&3.8922\\
		1600 &&& 9.37e-05 & 2.9880 &&& 6.13e-08 & 4.0029 &&& 7.69e-09 &  3.0114&&& 3.06e-08&3.9446\\
		3200 &&& 1.18e-05 & 2.9944 &&& 3.83 e-09 & 4.0001 &&& 9.58e-10 & 3.0057&&&1.96e-09 &3.9658\\
		6400 &&& 1.47e-06 & 2.9972 &&& 2.41e-10 & 3.9875 &&& 1.20e-10 & 3.0028&&& 1.20e-10&4.0245\\ 
		\hline \hline
	\end{tabular}
\end{center}
\label{Eg2rk3}
\end{table}
\begin{table}[h!]
\caption{Global error and order of convergence of RK4 and MQ-RBF-based RK4 methods applied to Example~\ref{6.2}.}
\begin{center}
	\begin{tabular}{ccccccccccccccccccc}
		\hline\hline
		&&& \multicolumn{2}{c}{\textbf{RK4 (C1)}}  &&& \multicolumn{2}{c}{\textbf{MQ-RBF--RK4 (C1$+$)}} &&& \multicolumn{2}{c}{\textbf{MQ-RBF--RK4 (C1$-$)}}& &&\\
		\cline{4-5}\cline{8-9}\cline{12-13}
		$N$ &&& error & order &&& error & order &&&  error & order &&&\\
		\hline\hline
		200  &&& 6.19e-04 & ... &&& 2.59e-05& ...   &&& 2.71e-04 &  ... \\
		400  &&& 3.98e-05 & 3.9582 &&& 8.30e-07 & 4.9696 &&& 9.21e-06 &  4.8785\\
		800  &&& 2.52e-06 & 3.9801 &&& 2.63e-08 & 4.9816 &&& 3.02e-07 &  4.9320\\
		1600 &&& 1.59e-07 & 3.9902 &&& 8.21e-10 & 5.0006 &&& 9.74e-09 &  4.9532\\
		3200 &&& 9.95e-09 & 3.9972 &&& 2.42e-11 & 5.0842 &&& 3.16e-10 & 4.9466\\
		6400 &&& 6.20e-10 & 4.0039 &&& 5.23e-12 & 2.2087 &&& 3.51e-12 & 6.4928\\ 
		\hline\hline
		&&& \multicolumn{2}{c}{\textbf{RK4 (C2)}}  &&& \multicolumn{2}{c}{\textbf{MQ-RBF--RK4 (C2$+$)}} &&& \multicolumn{2}{c}{\textbf{MQ-RBF--RK4 (C2$-$)}} \\
		\cline{4-5}\cline{8-9} \cline{12-13}
		$N$ &&& error & order &&& error & order &&&  error & order \\
		\hline\hline
		200  &&& 6.59e-04 & ... &&&2.59e-06 & ...   &&& 3.30e-03 &  ... \\
		400  &&& 4.27e-05 & 3.9485 &&& 7.98e-08 & 5.0204 &&& 1.14e-04 &  4.8526\\
		800  &&& 2.71e-06 & 3.9749 &&& 2.47e-09 & 5.0160 &&& 3.76e-06 &  4.9262\\
		1600 &&& 1.71e-07 & 3.9877 &&& 7.57e-11& 5.0273 &&& 1.21e-07 &  4.9595\\
		3200 &&& 1.07e-08 & 3.9966 &&& 1.72e-12 & 5.4619 &&& 3.83e-09 & 4.9790\\
		6400 &&& 6.71e-10 & 4.0008 &&& 6.63e-12 & -1.9491 &&& 1.09e-10 & 5.1374\\ 
		\hline\hline
	\end{tabular}
\end{center}
\label{Eg2rk4}
\end{table}

\end{example}

\begin{example}\label{6.3}
As another example, we consider the following stiff problem:
\begin{equation*}
\frac{du}{dt}=\frac{2t^2-u}{t^2u-t}, ~~~t\geq1,
\end{equation*}
with the initial condition $u(1)=2$. The exact solution for this problem is given by $u(t)=\frac{1}{t}+\sqrt{\frac{1}{t^2}+4t-4}$. We approximate the numerical solution at time $T=2$. 

\begin{table}[h!]
\caption{Optimal value of the shape parameter $\varepsilon_{2i}^2$ when the MQ-RBF--RK methods are applied to Example~\ref{6.3}.\footnote{The expressions for $\mathcal{N}_{32a}$, $\mathcal{D}_{32a}$, $\mathcal{N}_{32b}$, $\mathcal{D}_{32b}$, $\mathcal{N}_{34}$, $\mathcal{D}_{34}$, $\mathcal{N}_{41}$ and $\mathcal{N}_{42}$ are provided in Appendix.}}
\centering
\begin{tabular}{l cccc}
	\hline\hline
	\textbf{Methods}& &&
	$\varepsilon_{2i}^2$\\
	\hline
	MQ-RBF--RK2 \hspace{0.5cm}& &&$(-4t^4 + 2tu^3 - 3u^2 + 4t)/(t(tu - 1)^3u)$\\
	MQ-RBF--RK3 (\textbf{B1})  \hspace{0.5cm} &&&$-\frac{(4t^6 - 3t^2u^2 - 4t^3 + 6tu - 2)(4t^4 - 2tu^3 + 3u^2 - 4t)}{t(2t^4u^2 + 4t^5 - 3t^2u^3 - 4t^3u + 5tu^2 - 2t^2 - u)(tu - 1)^3}$\\
	MQ-RBF--RK3 (\textbf{B2 (a)}) &&& $\mathcal{N}_{32a}/ \mathcal{D}_{32a}$\\
	MQ-RBF--RK3 (\textbf{B2 (b)})  &&& $\mathcal{N}_{32b}/\mathcal{D}_{32b}$\\
	MQ-RBF--RK3 (\textbf{B3 (a)}) &&&$\frac{4 - 16t^9 + 32(u^3 + 1)t^6 + 12(-u^5 - 6u^2)t^5 + 12(3u^4 + 4u)t^4 + 8(-7u^3 - 4)t^3 + 48t^2u^2 + 3(u^4 - 4u)t - 4u^3}{(12t^6u + 2t^4u^2 + 4t^5 - 3t^2u^3 - 16t^3u + 5tu^2 - 2t^2 + 2u)(tu - 1)^3}$\\
	MQ-RBF--RK3 (\textbf{B3 (b)})  &&&$\frac{3(-16t^10 + 32t^7 + 4(u^5 + 2u^2)t^6 + 4(-3u^4 - 4u)t^5 + 8(u^3 - 2)t^4 + (-u^4 + 4u)t^2 + 4(u^3 + 1)t - 4u^2)}{(tu - 1)^4(4t^5 - 2t^3u + 3tu^2 - 2t^2 - 2u)t}$\\
	MQ-RBF--RK3 (\textbf{B4})  \hspace{0.5cm} &&&$\mathcal{N}_{34}/\mathcal{D}_{34}$\\
	MQ-RBF--RK4 (\textbf{C1$+$})  \hspace{0.5cm} &&&$ \mathcal{N}_{41}/(224(t^3 - 1/2)(tu - 1)^3t^3u^2) $\\
	MQ-RBF--RK4 (\textbf{C1$-$})  \hspace{0.5cm} &&&$conj(\mathcal{N}_{41})/(224(t^3 - 1/2)(tu - 1)^3t^3u^2)$\\
	MQ-RBF--RK4 (\textbf{C2$+$})  \hspace{0.5cm} &&&$ \mathcal{N}_{42}/(4(t^3 - 1/2)(tu - 1)^3t^3u^2) $\\
	MQ-RBF--RK4 (\textbf{C2$-$})  \hspace{0.5cm} 
	&&&$ conj(\mathcal{N}_{42})/(4(t^3 - 1/2)(tu - 1)^3t^3u^2)  $\\
	\hline\hline
\end{tabular}
\label{tabs3}
\end{table}
For all proposed methods, the values of $\varepsilon_{2i}^2$ are listed in Table~\ref{tabs3}. The computed rates of convergence for each MQ-RBF-based RK method, along with those of the corresponding classical RK methods for various values of $N$, are presented in Tables~\ref{Eg3rk2}, \ref{Eg1rk3}, and \ref{Eg3rk4}. These also include the global errors corresponding to each value of $N$.        
We can observe from these tables that the $s$-stage MQ-RBF RK methods are more accurate than the corresponding $s$-stage classical RK methods when using the same number of steps. However, some methods experience order reduction for specific values of $N$; for example, the MQ-RBF--RK3 (\textbf{B1}) method with $N=40$ and with $N=160$, the MQ-RBF--RK3 (\textbf{B3(a)}) method with $N=80$. On the other hand, some methods exhibit superconvergence, such as the MQ-RBF--RK3 (\textbf{B2(a)} method with $N=40$, anf the MQ-RBF--RK3 (\textbf{B3(a)} method with $N=40$. Moreover, except for the MQ-RBF-RK3 (\textbf{B1}) method with $N=20$ and $40$, and the MQ-RBF--RK4 (\textbf{C2$+$}) method with $N=20$, all the MQ-RBF--RK methods demonstrate lower global error than their corresponding RK counterparts. In addition to this, in terms of the order of convergence, most of the $s$-stage MQ-RBF--RK methods outperform even the classical RK methods with $s+1$ stages.

\begin{table}[h!]
\caption{Global error and order of convergence of RK2 and MQ-RBF-based RK2 methods applied to Example~\ref{6.3}.}
\begin{center}
	\begin{tabular}{ccccccccccccccc}
		\hline\hline
		&&& \multicolumn{2}{c}{\textbf{RK2}}  &&& \multicolumn{2}{c}{\textbf{MQ-RBF--RK2}} &&& \\
		\cline{3-5}\cline{7-9} 
		$N$    && &     error    &    order &&&   error  &    order   \\
		\hline\hline			
		20  && & 1.56e-04  & .... &&& 2.03e-05 & .... \\
		40  & && 3.80e-05  & 2.0399 &&& 2.44e-06 & 3.0552 \\
		80 & && 9.38e-06   & 2.0194 &&& 2.99e-07 & 3.0280\\
		160 & && 2.33e-06  & 2.0096 &&& 3.71e-08& 3.0141\\
		320 & && 5.80e-07  & 2.0048 &&&4.61e-09 & 3.0071\\
		\hline\hline
	\end{tabular}
\end{center}
\label{Eg3rk2}
\end{table}
\begin{table}[h!]
\caption{Global error and order of convergence of RK3 and MQ-RBF-based RK3 methods applied to Example~\ref{6.3}.}
\begin{center}
	\begin{tabular}{cccccccccccccccccccccc}
		\hline\hline
		&&& \multicolumn{2}{c}{\textbf{RK3 (B1)}}  &&& \multicolumn{2}{c}{\textbf{MQ-RBF--RK3 (B1)}} &&& \multicolumn{2}{c}{\textbf{RK3 (B2(a))}}& &&\multicolumn{2}{c}{\textbf{MQ-RBF--RK3 (B2(a))}}\\
		\cline{4-5}\cline{8-9}\cline{12-13}\cline{16-17}
		$N$ &&& error & order &&& error & order &&&  error & order &&&error & order\\
		\hline\hline
		20  &&& 7.68e-07 & ... &&& 9.98e-07 & ...   &&& 4.40e-06  &  ... &&& 4.21e-05&...\\
		40  &&& 1.18e-07 & 2.6964 &&&  1.95e-07 &  2.3535 &&&5.36e-07 &  3.0362&&&1.27e-06 &5.0466\\
		80  &&& 1.60e-08 & 2.8921 &&& 1.13e-08  & 4.1138  &&& 6.61e-08 &  3.0200&&&4.23e-08 &4.9110\\
		160 &&& 2.06e-09 & 2.9545 &&& 7.64e-08  &  -2.7606 &&& 8.20e-09 &  3.0105&&& 1.64e-09&4.6873\\
		320 &&& 2.61e-10 & 2.9791 &&&  2.42e-09 &  4.9811 &&& 1.02e-09 & 3.0054&&& 7.62e-11&4.4310\\
		\hline\hline
		&&& \multicolumn{2}{c}{\textbf{RK3 (B2(b))}}  &&& \multicolumn{2}{c}{\textbf{MQ-RBF--RK3 (B2(b))}} &&& \multicolumn{2}{c}{\textbf{RK3 (B3(a))}}& &&\multicolumn{2}{c}{\textbf{MQ-RBF--RK3 (B3(a))}}\\
		\cline{4-5}\cline{8-9}\cline{12-13}\cline{16-17}
		$N$ &&& error & order &&& error & Order &&&  error & order &&&error & order\\
		\hline\hline
		20  &&& 3.80e-06 & ... &&& 2.33e-07& ...   &&& 7.04e-06 &  ... &&& 9.38e-07&...\\
		40  &&& 4.66e-07 & 3.0259 &&& 1.37e-08 & 4.0857 &&& 8.83e-07 &  2.9948&&& 1.13e-09&9.7025\\
		80  &&& 5.76e-08 & 3.0154 &&& 8.32e-10 & 4.0412 &&& 1.11e-07 &  2.9983&&& 1.51e-09&-0.4283\\
		160 &&& 7.16e-09& 3.0083 &&& 5.13e-11 & 4.0201 &&& 1.38e-08 &  2.9994&&&1.32e-10 &3.5233\\
		320 &&& 8.93e-10 & 3.0043 &&& 3.18e-12 & 4.0097 &&& 1.73e-09 & 2.9997&&& 9.25e-12&3.8314\\ 
		\hline\hline
		&&& \multicolumn{2}{c}{\textbf{RK3 (B3(b))}}  &&& \multicolumn{2}{c}{\textbf{MQ-RBF--RK3 (B3(b))}} &&& \multicolumn{2}{c}{\textbf{RK3 (B4)}}& &&\multicolumn{2}{c}{\textbf{MQ-RBF--RK3 (B4)}}\\
		\cline{4-5}\cline{8-9}\cline{12-13}\cline{16-17}
		$N$ &&& error & Order &&& error & order &&&  error & order &&&error & order\\
		\hline\hline
		20  &&& 4.99e-06 & ... &&& 2.89e-07& ...   &&& 6.78e-06 &  ... &&& 9.43e-07&...\\
		40  &&& 6.08e-07 & 3.0372 &&& 1.74e-08 & 4.0554 &&& 8.36e-07 &  3.0185&&&5.55e-08 &4.0857\\
		80  &&& 7.50e-08 & 3.0201 &&& 1.07e-09 & 4.0257 &&& 1.04e-07 &  3.0101&&& 3.37e-09&4.0428\\
		160 &&& 9.31e-09 & 3.0104 &&& 6.62e-11 & 4.0124 &&& 1.30e-08 &  3.0052&&& 2.07e-10&4.0214\\
		320 &&& 1.16e-09 & 3.0053 &&& 4.12e-12 & 4.0064 &&& 1.61e-09 & 3.0027&&&1.29e-11 &4.0107\\
		\hline \hline
	\end{tabular}
\end{center}
\label{Eg3rk3}
\end{table}
\begin{table}[h!]
\caption{Global error and order of convergence of RK4 and MQ-RBF-based RK4 methods applied to Example~\ref{6.3}.}
\begin{center}
	\begin{tabular}{ccccccccccccccccccc}
		\hline\hline
		&&& \multicolumn{2}{c}{\textbf{RK4 (C1)}}  &&& \multicolumn{2}{c}{\textbf{MQ-RBF--RK4 (C1$+$)}} &&& \multicolumn{2}{c}{\textbf{MQ-RBF--RK4 (C1$-$)}}& &&\\
		\cline{4-5}\cline{8-9}\cline{12-13}
		$N$ &&& error & order &&& error & order &&&  error & order &&&\\
		\hline\hline
		20  &&& 1.51e-08 & ... &&& 7.88e-09& ...   &&& 7.04e-09 &  ... \\
		40  &&& 1.22e-09 & 3.6299 &&& 2.45e-10 & 5.0089 &&& 2.21e-10 &  4.9914\\
		80  &&& 8.46e-11 & 3.8520 &&& 7.62e-12 & 5.0056 &&& 6.93e-12 &  4.9980\\
		160 &&& 5.54e-12 & 3.9328 &&& 2.37e-13 & 5.0052 &&& 2.16e-13 &  5.0016\\
		320 &&& 3.52e-13 & 3.9759 &&& 8.44e-15 & 4.8128 &&& 7.99e-15 & 4.7579\\
		\hline\hline
		&&& \multicolumn{2}{c}{\textbf{RK4 (C2)}}  &&& \multicolumn{2}{c}{\textbf{MQ-RBF--RK4 (C2$+$)}} &&& \multicolumn{2}{c}{\textbf{MQ-RBF--RK4 (C2$-$)}} \\
		\cline{4-5}\cline{8-9} \cline{12-13}
		$N$ &&& error & order &&& error & order &&&  error & order \\
		\hline\hline
		20  &&& 9.04e-09 & ... &&&1.18e-08 & ...   &&& 6.71e-09 &  ... \\
		40  &&& 5.56e-10 & 4.0234 &&& 1.25e-10 & 5.1790 &&& 1.67e-10 &  5.3282\\
		80  &&& 3.44e-11 & 4.0152 &&& 9.54e-12 & 5.0888 &&& 4.77e-12 &  5.1315\\
		160 &&& 2.14e-12 & 4.0080 &&& 2.90e-13& 5.0401 &&& 1.40e-13 &  5.0864\\
		320 &&& 1.35e-13 & 3.9845 &&& 8.45e-15 & 5.1018 &&& 4.00e-15 & 5.1339\\
		\hline\hline
	\end{tabular}
\end{center}
\label{Eg3rk4}
\end{table}

\end{example}
\begin{example}\label{6.4}
Next, we consider a system of IVPs to demonstrate that the proposed MQ-RBF-based RK2 method can be effectively applied to systems of ODEs, exhibiting the same improvements in accuracy and order as observed for scalar equations. Consider the following system:
\begin{eqnarray*}
\frac{du}{dt} =\begin{pmatrix}
	e^t\\0
\end{pmatrix}-\begin{pmatrix}
	5 &-3\\
	3&-1
\end{pmatrix} u,
\end{eqnarray*}
where $u=\begin{pmatrix}
u_1\\u_2   \end{pmatrix}$, and the  initial condition is given by $u(0)=\begin{pmatrix}
1\\0  \end{pmatrix} $. The exact solution of this system is $u(t)= \begin{pmatrix}
(1-2t)e^{-2t}\\ \left(\frac{1}{3}-2t\right)e^{-2t}-\frac{e^t}{3}
\end{pmatrix}$. We compute the numerical solution of this system at the final time $T= 5$, and compare it with the exact solution at the same point to determine the error.

\begin{table}[h!]
\caption{Global error and order of convergence of RK2 and MQ-RBF-based RK2 methods applied to Example~\ref{6.4}.}
\begin{center}
	\begin{tabular}{ccccccccccccccc}
		\hline\hline
		&&& \multicolumn{2}{c}{\textbf{RK2}}  &&& \multicolumn{2}{c}{\textbf{MQ-RBF--RK2}} &&& \\
		\cline{3-5}\cline{7-9} 
		$N$    && &     error    &    order &&&   error  &    order   \\
		\hline\hline			
		20  && & 3.87e-01  & .... &&& 4.39e-02 & .... \\
		40  & && 7.17e-02  & 2.4309 &&& 3.93e-03 & 3.4798 \\
		80 & && 1.62e-02   & 2.1482 &&& 4.53e-04 & 3.1193\\
		160 & && 3.90e-03  & 2.0528 &&& 5.55e-05& 3.0294\\
		320 & && 9.61e-04  & 2.0211 &&&6.91e-06 & 3.0042\\
		\hline\hline
	\end{tabular}
\end{center}
\label{Eg4rk2}
\end{table}

The global error at the final time and the corresponding order of accuracy for different values of $N$ are shown in Table~\ref{Eg4rk2} for both the standard and MQ-RBF-based RK2 methods. Note that, in computing the shape parameter for this system, the term $f_u$ refers to the Jacobian of the vector field $f$ with respect to $u$, rather than the ordinary derivative used in the scalar case. Care must be taken with the order of terms in the product involving $f$ and $f_u$, to ensure that the resulting expression is mathematically well-defined. Furthermore, unlike the scalar case where the error is computed as the absolute difference between the numerical and exact solutions, the error for systems of equations is defined as the norm of the component-wise difference between the numerical solution and the exact solution vectors. It is evident from Table~\ref{Eg4rk2} that the classical RK2 method exhibits only second-order accuracy, whereas the MQ-RBF-based RK2 method achieves third-order accuracy. It is also observed that, for the same number of steps $N$, the MQ-RBF-based RK2 method yields a smaller global error than the standard RK2 method.
\end{example}
\begin{example}\label{6.5}
As a final example, we consider the Duffing equation, a second-order nonlinear differential equation that describes certain damped and driven oscillators, given as
\begin{equation}\label{fe}
q''+\omega^2q=k^2(2q^3-q), \ \ q(0)=0, \ \ q'(0)=\omega, \quad 0<t\le 20,  
\end{equation}
where $0\le k<\omega$. Here, $\omega$ represents the natural frequency of the system, $k$ is the nonlinearity parameter, $q(t)$ denotes the displacement over time $t$, and the right-hand side of the equation introduces a nonlinear restoring force. The analytic solution of \eqref{fe} is given as 
$q(t)=sn(\omega t,k/\omega)$,
where $sn$ is the Jacobi elliptic sine function with modulus $(k/\omega)^2$. This function describes a periodic and bounded motion influenced by nonlinearity.

To apply numerical methods to IVPs, we convert this second-order ODE to a system of first-order ODEs. Define $u=\begin{pmatrix}
p\\q
\end{pmatrix}$, where $p=q'$, then \eqref{fe} get reduced to
$$u'=f(u):=\begin{pmatrix}
-\omega^2q+k^2(2q^3-q)\\p
\end{pmatrix}, \ \ \text{with} \ \ u(0)=\begin{pmatrix}
\omega \\0
\end{pmatrix}.$$
We solve the above system numerically over the interval $[0,20]$ using both the classical RK2 method and the MQ-RBF-based RK2 method, with parameters $k=0.03$, and $\omega=10$. For different values of $N$, the global error and order of convergence of the RK2 and MQ-RBF-RK2 methods at the final time $t=20$ are presented in Table~\ref{Eg5rk2}. It is evident that the MQ-RBF-based RK2 method outperforms the classical RK2 method. Although both methods require small step sizes for the numerical solution to closely approximate the exact solution, the MQ-RBF-based RK2 method consistently yields lower errors. As illustrated in Figure~\ref{fig2}, for a step size $h=0.06$, the energy obtained using the classical RK2 and MQ-RBF-based RK2 methods increases over time, leading to potential blow-up during long-term integration. However, the rate of energy growth is significantly lower for the MQ-RBF-based RK2 method, making it closer to the exact solution over longer time intervals compared to the classical RK2 method. 
\begin{table}[h!]
\caption{Global error and order of convergence of RK2 and MQ-RBF-based RK2 methods applied to Example~\ref{6.5}.}
\begin{center}
	\begin{tabular}{ccccccccccccccc}
		\hline\hline
		&&& \multicolumn{2}{c}{\textbf{RK2}}  &&& \multicolumn{2}{c}{\textbf{MQ-RBF--RK2}} &&& \\
		\cline{3-5}\cline{7-9} 
		$N$    && &     error    &    order &&&   error  &    order   \\
		\hline\hline			
		640  && & 2.73e-00  & .... &&& 5.20e-01 & .... \\
		1280  & && 5.96e-01  & 2.1955 &&& 5.15e-02 & 3.3356 \\
		2560 & && 1.07e-01   & 2.4780 &&& 6.07e-03 & 3.0860\\
		5120 & && 2.46e-02  & 2.1213 &&& 7.41e-04& 3.0347\\
		10240 & && 6.61e-03  & 2.0123 &&&9.15e-05 & 3.0173\\		
		\hline\hline
	\end{tabular}
\end{center}
\label{Eg5rk2}
\end{table}
\begin{figure}[h!]
\centering
\includegraphics[width=0.45\linewidth]{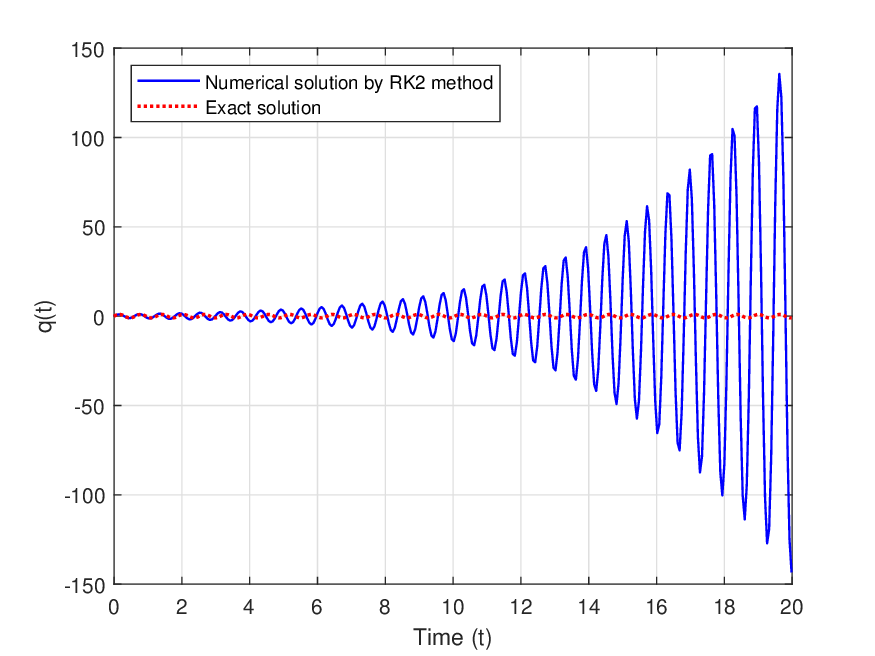}
\includegraphics[width=0.45\linewidth]{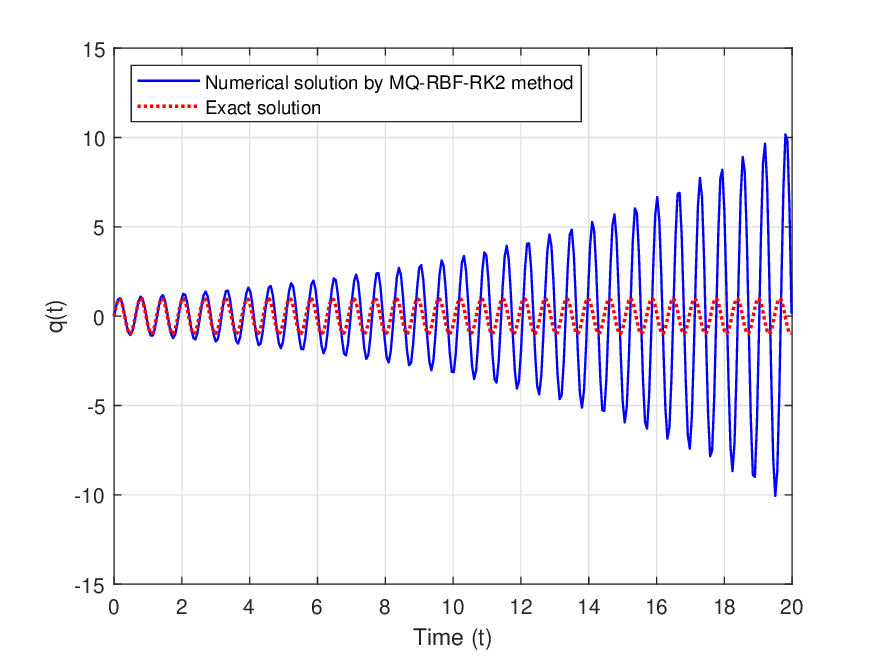}
\caption{Solution profile of the Duffing equation on the interval $[0,20]$ using the RK2 method (left) and the MQ-RBF-based RK2 method (right), with a step size of $h=0.06$.}
\label{fig2} 
\end{figure}
\end{example}
\section{Conclusion}\label{sec6}
\noindent
In this article, we have developed multiquadratic radial basis function (MQ-RBF)–based Runge-Kutta (RK) methods to solve initial value problems (IVPs). By incorporating MQ-RBFs with appropriately chosen shape parameters into the classical explicit RK methods, we have enhanced their accuracy and order. While a classical $s$-stage RK method typically achieves an order of accuracy $s$, the proposed $s$-stage MQ-RBF-based RK method analytically achieves an improved order of accuracy $s+1$. We have provided theoretical results establishing the convergence of the proposed methods. Furthermore, we have plotted and analyzed the stability regions of the MQ-RBF-based RK methods and compared them with those of the classical RK methods. Numerical experiments demonstrate the superior performance and potential of the proposed methods. We have considered both stiff and non-stiff test problems with varying numbers of time steps $N$. For non-stiff problems, the $s$-stage MQ-RBF-based RK methods consistently achieve the $s+1$-order convergence. However, for non-stiff problems, some MQ-RBF-based RK methods exhibited order reduction for certain values of $N$. Nevertheless,  for both stiff and non-stiff problems, the global errors produced by the MQ-RBF-based RK methods were generally lower than those of their RK counterparts. In many cases, the proposed $s$-stage MQ-RBF-based RK methods even outperform the classical $s+1$-stage RK methods in terms of order or accuracy, or both. Therefore, we conclude that MQ-RBF-based RK methods are more suitable for solving IVPs compared to classical methods. We have also considered examples, including the Duffing equation, to demonstrate that the MQ-RBF-based RK2 method can be readily extended to systems of IVPs. In future work, one may consider developing implicit or semi-implicit variants of the proposed methods to enhance their performance for stiff problems. Moreover, generalizations involving arbitrarily chosen shape parameters and extensions to non-uniform meshes present promising directions for further research. The extension of MQ-RBF-based RK methods to higher-dimensional problems also represents a worthwhile area for future exploration.

\section*{ Acknowledgments } 
\noindent
The work of \textbf{Rajesh Yadav} is supported by the Council of Scientific and Industrial Research, India, with file no. 09/1201(13064)/2022-EMR-I. The work of \textbf{Alpesh Kumar} is supported by the Science and Engineering Research Board, India, under the MATRICS scheme with sanction order number MTR/2022/000149.

\section*{Appendix}
\begin{eqnarray*}
\mathcal{N}_{32a}&:=& -11067267806781600t^{10} + (10456385037949128u^3 + 22134535613563200)t^7 + (-2461375567279164u^5\\
&&- 18145953124202856u^2)t^6 + 7384126701837492 u(u^3 + 2/3)t^5 +(-15379136172507456u^3 - 13834084758477000)t^4\\
&&+ (-1230687783639582u^5 + 13223201989644528u^2)t^3 + 4922751134558328u(u^3 + 1/2)t^2 + (-4770030442350210u^3\\
&&+ 305441384416236)t +  386262853507614u^2\\ 
\mathcal{D}_{32a}&:=&t(3587275474121788t^6u + 820458522426388t^4u^2+ 1640917044852776t^5 - 1230687783639582t^2u^3 - \\
&&5228192518974564t^3u + 2051146306065970tu^2 - 820458522426388t^2 + 486589607317253u)(tu - 1)^3\\
\mathcal{N}_{32b}&:=&-2443531075329888t^{10} + (-323285876365512u^3 + 4887062150659776)t^7 + (772525707015228u^5\\
&&+ 1257454521563496u^2)t^6 - 2317577121045684u(u^3 + 2/3)t^5 + (1868337290395968u^3 - 3054413844162360)t^4\\
&&+ (386262853507614u^5 + 287596892466960u^2)t^3 - 1545051414030456u(u^3 + 1/2)t^2 + (2236755651954306u^3\\
&&+ 1383408475847700)t - 1230687783639582u^2\\
\mathcal{D}_{32b}&:=&t(353374199827396t^6u - 257508569005076t^4u^2 - 515017138010152t^5 + 386262853507614t^2u^3 + 161642938182756t^3u\\
&&- 643771422512690tu^2 + 257508569005076t^2 + 217097834459387u)(tu - 1)^3\\
\mathcal{N}_{34}&:=&-96t^{10} + (16u^3 + 192)t^7 + (-4u^5 + 24u^2)t^6 + (8u^7 + 44u^4 - 96u)t^5 + (-36u^6 - 152u^3 - 88)t^4 + (56u^5 + 160u^2)t^3 \\
&&+ (-32u^4 - 20u)t^2 + (12u^3 + 20)t - 14u^2\\
\mathcal{D}_{34}&:=&t(12t^6u - 2t^4u^2 - 4t^5 + 3t^2u^3 - 8t^3u - 5tu^2 + 2t^2 + 4u)(tu - 1)^3\\
\mathcal{N}_{41}&:=&sqrt\Big((98560t^{12} - 49280t^9)u^{11} + (-591360t^{11} + 313104t^8)u^{10} + (-197120t^{13} + 1486496t^{10} - 815936t^7)u^9 \\
&&+ (1250960t^{12} - 2238848t^9 + 1159520t^6)u^8 + (-78848t^{14} - 3514368t^{11} + 3065040t^8 - 1166792t^5)u^7 + (652608t^{13} \\
&&+ 5110192t^{10} - 4328488t^7 + 1171024t^4)u^6 + (452864t^{15} - 2483392t^{12} - 3147232t^9 + 4182480t^6 - 1022312t^3)u^5 \\
&&+ (-547072t^{14} + 2472576t^{11} + 890868t^8 - 2387020t^5 + 512961t^2)u^4 + (-140800t^{16} - 49152t^{13}+ 420736t^{10} \\
&&- 1702992t^7 + 1213664t^4 - 121660t)u^3 + (-793600t^{18} + 2630016t^{15} - 2943776t^{12} + 757952t^9 + 973880t^6 \\
&&- 472232t^3 + 12100)u^2 + (294400t^{17} - 959232t^{14} + 1110976t^{11} - 383168t^8 - 172272t^5 + 58960t^2)u \\
&&+ 33856(t^6 - 91/46t^3 + 67/46)^2t^4\Big) - 800t^9u - 184t^8 + 88t^7u^2 + (268u^4 + 1024u)t^6 + (-24u^3 + 364)t^5 \\
&&+ (-132u^5 - 886u^2)t^4 + (328u^4 + 172u)t^3 + (-516u^3 - 268)t^2 + 553tu^2 - 110u\\
\mathcal{N}_{42}&:=&sqrt\Big((160t^{12} - 80t^9)u^{11} + (-960t^{11} + 624t^8)u^{10} + (-320t^{13} + 1904t^{10} - 1880t^7)u^9 + (2240t^{12} - 1568t^9 + 3140t^6)u^8\\
&&+ (448t^{14} - 7392t^{11} + 3384t^8 - 4220t^5)u^7 + (9988t^{10} - 9988t^7 + 5521t^4)u^6 + (-256t^{15} - 3520t^{12} - 1216t^9 + 11880t^6 \\
&&- 5180t^3)u^5 + (4064t^{14} - 1344t^{11} - 2472t^8 - 7336t^5 + 2688t^2)u^4 + (-640*t^16 - 4224t^{13} + 7552t^{10} - 4944t^7 + 4400t^4\\
&&- 760t)u^3 + (-1984t^{18} + 7488t^{15} - 11120t^{12} + 6176t^9 + 1424t^6 - 1472t^3 + 100)u^2 + (-1280t^{17} + 2688t^{14} + 256t^{11} \\
&&- 2048t^8 - 192t^5 + 160t^2)u + 256(t^6 - t^3 - 1/2)^2t^4\Big) - 40t^9u + 16t^8 + 8t^7u^2 + (8u^4 + 44u)t^6 + (6u^3 - 16)t^5 \\
&&+ (-12u^5 - 56u^2)t^4 + (38u^4 + 32u)t^3 + (-51u^3 - 8)t^2 + 38tu^2 - 10u.
\end{eqnarray*}


\end{document}